\documentclass[preprint,11pt,authoryear]{elsarticle}

\usepackage{amssymb}

\usepackage{setspace}
\usepackage{amsfonts,amssymb} 
\usepackage{mathrsfs} 
\usepackage{lipsum}
\usepackage{graphicx}
\usepackage{bm}
\usepackage{amsfonts, amssymb, amsmath}
\usepackage{multicol}
\usepackage{multirow}
\usepackage{lineno,hyperref}
\usepackage{autobreak}
\usepackage{color}
\usepackage{epstopdf}
\usepackage{mathrsfs}
\ifpdf
\DeclareGraphicsExtensions{.eps,.pdf,.png,.jpg}
\else
\DeclareGraphicsExtensions{.eps}
\fi
\newproof{proof}{Proof}
\newtheorem{thm}{Theorem}
\newtheorem{cor}{Corollary}

\newtheorem{lem}{Lemma}

\newtheorem{asm}{Assumption}
\newtheorem{defn}{Definition}

\newcommand{\eq}{\begin{equation}}
	\newcommand{\nq}{\end{equation}}
\newcommand{\eqa}{\begin{eqnarray}}
	\newcommand{\nqa}{\end{eqnarray}}

\newcommand{\nn}{\nonumber}
\newcommand{\hb}{\hfill $\Box$}

\usepackage{algorithm}
\usepackage{algorithmic}

\makeatletter
\newenvironment{breakablealgorithm}
{
		\begin{center}
			\refstepcounter{algorithm}
			\hrule height.8pt depth0pt \kern2pt
			\renewcommand{\caption}[2][\relax]{
				{\raggedright\textbf{\ALG@name~\thealgorithm} ##2\par}%
				\ifx\relax##1\relax 
				\addcontentsline{loa}{algorithm}{\protect\numberline{\thealgorithm}##2}%
				\else 
				\addcontentsline{loa}{algorithm}{\protect\numberline{\thealgorithm}##1}%
				\fi
				\kern2pt\hrule\kern2pt
			}
		}{
		\kern2pt\hrule\relax
	\end{center}
}
\makeatother

\journal{European Journal of Operational Research}

\begin{document}
	
	\begin{frontmatter}
		
		
		
    \title{Restricted inverse optimal value problem on linear programming under weighted $l_1$ norm}
		\tnotetext[t1]{Research supported by the National Natural  Science Foundation of China (11471030).}
		
		\author[1]{Junhua Jia}
		\ead{230218184@seu.edu.cn}
		
		\author[1]{Xiucui Guan\corref{cor1}}

		\author[1]{Xinqiang Qian}
  \ead{230208660@seu.edu.cn}
		\author[2]{Panos M. Pardalos}
   \ead{pardalos@ufl.edu}
		\cortext[cor1]{Corresponding author: xcguan@163.com}
 

		\affiliation[1]{organization={School of Mathematics, Southeast University},
			addressline={No.2 Sipailou},
			postcode={210096},
			city={Nanjing},
			country={China}}

		\affiliation[2]{organization={Center for Applied Optimization, University of Florida},
			addressline={P.O. Box 116595},
			postcode={32611},
			city={Gainesville, Florida},
			country={USA}}
		
		
		\begin{abstract}
			We study the restricted inverse optimal value problem on linear programming under weighted $l_1$ norm (\textit{RIOVLP}$_1$). Given a linear programming problem \textit{LP}$_c: \min \{cx|Ax=b,x\geq 0\}$ with a feasible solution $x^0$ and a value $K$, we aim to adjust the vector $c$ to $\bar{c}$ such that $x^0$ becomes an optimal solution of the problem \textit{LP}$_{\bar c}$ whose objective value $\bar{c}x^0$ equals $K$. The objective is to minimize the distance $\|\bar c - c\|_1=\sum_{j=1}^nd_j|\bar c_j-c_j|$  under weighted $l_1$ norm. 
			Firstly, we formulate the problem (\textit{RIOVLP}$_1$) as a linear programming problem by dual theories. 
               Secondly, we construct a sub-problem $(D^z)$, which has the same form as $LP_c$,  of the dual  (\textit{RIOVLP}$_1$) problem corresponding to a given value $z$. 
                Thirdly, when the coefficient matrix $A$ is unimodular, we design a binary search algorithm to calculate the critical value $z^*$ corresponding to an optimal solution of the problem (\textit{RIOVLP}$_1$).  
               Finally, we solve the (\textit{RIOV}) problems on Hitchcock and shortest path problem, respectively, in $O(T_{MCF}\log\max\{d_{max},x^0_{max},n\})$ time, where  we solve a sub-problem $(D^z)$ by minimum cost flow in $T_{MCF}$ time in each iteration. The values $d_{max},x^0_{max}$ are the maximum values of $d$ and $x^0$, respectively.
		\end{abstract}
		
		
		
		
	\end{frontmatter}
	
	
	\section{Introduction} \label{introduction}
	Since \cite{Burton_1992} first studied the inverse shortest path problem, many researchers have considered different inverse combinatorial optimization problems such as inverse spanning tree problem, inverse shortest path problem, inverse minimum cost flow problem, inverse minimum cut problem, inverse maximum matching problem. The inverse combinatorial optimization problems have broad applications which can be found in \cite{Mohajerin_2018}, \cite{Heuberger_2004}, \cite{Ahuja_Orlin_2001} and \cite{Burton_1992}.

 Let \textit{LP}$_{c}$ be a standard linear programming (\textit{LP})  problem,
	\eqa  
	& \min & cx\nn\\
	(\textit{LP}_c) \ & {\hbox {s.t.}}
	&  Ax=b, \label{eq-LP}\\
	& & x\geq 0.\nn
	\nqa
	where $A$ is an $m \times n$ matrix and $m<n$, $c^T$ and $x$ are $n \times 1$ vectors and $b$ is an $m \times 1$ vector. 
 
	As some combinatorial optimization problems can be described as \textit{LP} problems,  \cite{Zhang_Liu_1996} first studied the inverse \textit{LP} problem under unit $l_1$ norm (\textit{ILP}$_{u1}$). 
	Let $x^0$ be a given feasible solution of the problem (\textit{LP}). The aim of problem (\textit{ILP}$_{u1}$) is to minimize the modification $\|\bar c - c\|$ under unit $l_1$ norm such that $x^0$ becomes an optimal solution of the modified problem (\textit{LP}$_{\bar c}$). 
	They transformed the problem (\textit{ILP}$_{u1}$) into another \textit{LP} problem and extended their results to the inverse bounded \textit{LP} problem with a bounded variable constraint $ l\leq x\leq u$. Furthermore, they applied their research methods to the inverse minimum cost flow problem and inverse assignment problem under unit $l_1$ norm.
	\cite{Zhang_Liu_1999} continued to consider a special case of problem (\textit{ILP}$_{u1}$) in which the given feasible solution $x^0$ and one optimal solution of the original \textit{LP} problem are 0-1 vectors. They gave a method based on dual theories for solving this special case and applied the method to the inverse shortest path problem and inverse assignment problem under unit $l_1$ norm.
	\cite{Huang_Liu_1999} also studied the inverse bounded \textit{LP} problem and applied their results to the inverse minimum
	weight perfect $k$-matching problem on bipartite graphs under unit $l_1$ norm. 
\cite{Ahuja_Orlin_2001} studied the inverse canonical \textit{LP} problem   (\textit{ICLP}$_{1}$) under weighted $l_1$ norm.
They transformed the  problem (\textit{ICLP}$_{1}$) into an \textit{LP} problem and transformed the inverse shortest path problem, inverse minimum cut problem, inverse minimum cost flow problem and inverse assignment problem under weighted $l_1$ norm into some minimum cost flow problems. 
	\cite{Chan_kaw_2020} concentrated on imputing  unspecified constraint coefficient matrix $A$ and a cost vector for a given linear optimization problem. 
	\cite{Ghobadi_Mahmoudzadeh_2021} inferred the feasible region of \textit{LP} problem that would render the given solutions  feasible while making some optimal
	for the given cost function.  
	
	\cite{Ahmed_2005} studied the inverse optimal value problem (\textit{IOVLP}) on \textit{LP}. Given a desired optimal objective value $K$, and a set $C$ of feasible cost vectors in an (\textit{LP}), determine a cost vector $\bar{c}\in C$ such that the optimal objective value of the new problem \textit{LP}$_{\bar c}$ is closest to the desired value $K$. They proved the problem (\textit{IOVLP}) is NP-hard.  \cite{LV_2008} and \cite{LV_2010} studied this problem under more general conditions using a nonlinear bilevel programming approach. 
	
	In this paper, we will study the  restricted inverse optimal value problem (\textit{RIOVLP}$_1$) on \textit{LP} under weighted $l_1$ norm. Similar to the classical (\textit{ILP}) problem, its objective is to minimize the modification  $\|\bar c - c\|$ under weighted $l_1$ norm. Different to the classical (\textit{ILP}) problem, in (\textit{RIOVLP}$_1$) we not only require that the given feasible solution $x^0$ becomes an optimal solution of the problem (\textit{LP}$_{\bar c}$) but also require that the optimal objective value $\bar{c}x^0$ equals the given value $K$.  
	There are two main differences compared the problem (\textit{RIOVLP}$_1$) with (\textit{IOVLP}). One difference is on the optimization objectives. The problem (\textit{RIOVLP}$_1$) aims to minimize the distance $\|\bar c-c\|$, while the problem (\textit{IOVLP}) tries to minimize $|\bar c x^*-K|$ among $\bar{c}\in C$. The other difference is on the constraint conditions. In (\textit{RIOVLP}$_1$) we impose a constraint on the optimal value $\bar{c}x^0$, which is equal to the given value $K$, while in  (\textit{IOVLP}), there is no cadidate solution for consideration.
	
	Some restricted inverse optimal value problems on combinatorial optimization structures have been studied. \cite{Jia_Guan_2023}, \cite{Zhang_Guan_2021}, \cite{WangHui_2021} and \cite{Zhang_Guan_2020} considered the restricted inverse optimal value problems on minimum spanning tree under different norms and 
	proposed combinatorial algorithms in polynomial time.  \cite{ZhangQ_2023} studied the restricted inverse optimal value problem of shortest path on trees and devised an $O(n^2)$ algorithm under weighted $l_1$ norm and an $O(n)$ algorithm under unit $l_1$ norm. \cite{Zhang_Cai_1998} considered a more general restricted inverse optimal value problem under weighted $l_1$ norm on minimum cut which requires a set of cuts (not a cut) to become minimum cuts to make their objective value within a certain range ( not equal to a given value). \cite{Cui_2010} showed the restricted inverse optimal value problem on shortest path for general graphs is NP-hard when a collection of source-sink pairs with prescribed distances is given. 
	
This paper is organized as follows. In section \ref{section_RIOVLP_Properties}, we first formulate the problem (\textit{RIOVLP}$_1$) as an \textit{LP} problem by the dual theories.  Then we analyze some properties of a sub-problem $(D^z)$ of the dual   (\textit{RIOVLP}$_1$)  problem with respect to a given value $z$.
    In section \ref{section_RIOVLP_unimodular}, we design a binary search algorithm to calculate the critical value $z^*$ corresponding to an optimal solution of the problem (\textit{RIOVLP}$_1$) by solving a sub-problem in each iteration.
	In section \ref{section_applications}, we apply the above methods to the restricted inverse optimal value problems on Hitchcock and shortest path problem, respectively. Finally, we give some conclusions and further research in section \ref{section_conclusions}.
\setcounter{equation}{0}	
	\section{Properties of the problem (\textit{RIOVLP}$_1$) } \label{section_RIOVLP_Properties}
	In this section, we study the restricted inverse optimal value problem on linear programming under weighted $l_1$ norm. We first formulated the problem as an \textit{LP} problem, then analyze some properties of its sub-problems.
	
	\subsection{The mathematical model of the  problem (\textit{RIOVLP}$_1$)}
	Let $x^0$ be a given feasible solution, $\mathcal{F}_0=\{x\in R^n|Ax=b,x\geq 0\}$ be the feasible region of the problem (\textit{LP}) and $K$ be a given real number. We aim to adjust the vector $c$ to $\bar{c}$ such that $x^0$ becomes an optimal solution under $\bar{c}$ whose objective value $\bar{c}x^0$ equals $K$.  
Given  a $1\times n$ positive vector $d> 0$, the problem (\textit{RIOVLP}$_1$) under weighted $l_1$ norm can be formulated as follows.
	\eqa
	&\min &\sum_{j=1}^n d_j|\bar c_j-c_j|\nn\\
	(\textit{RIOVLP}_1) \ & {\hbox {s.t.}}  & \min_{x\in \mathcal{F}_0}\bar{c}x=K, \label{RIOVLP_w-1}\\
	&& \bar{c}x^0=K.\label{RIOVLP_w-2}
	\nqa
	 
	It follows from the constraint (\ref{RIOVLP_w-1}) that the problem (\textit{RIOVLP}$_1$) is not an \textit{LP} problem. Fortunately, we can turn it into an \textit{LP} problem by dual theories of \textit{LP}. Next, we will explain the process in details.
	
	Let us associate a dual variable $\pi\in R^m$ with the constraint (\ref{eq-LP}). Then the dual problem of (\textit{LP}$_{\bar c}$) can be stated as follows.
	\eqa  
	&  \max & \pi b\nn\\
	(\textit{DLP}_{\bar{c}}) \ & {\hbox {s.t.}}
	&  \pi A \leq \bar{c}. \label{DLP_bar_c-con-1} 
	\nqa
	
	Let $J=\{j|x^0_j=0\}$, $\bar{J}=\{j|x^0_j>0\}$ and $A_j$ be the $j$-th column of $A$. 
	
	\begin{thm} \label{model-linearization}
		If $(\pi^*,\bar c^*)$ is an optimal solution of the problem below, 
		\eqa
		\min &&\sum_{j=1}^n d_j|\bar c_j-c_j| \nn\\
		(\textit{RIOVLP}_{1}^1) \ {\hbox {s.t.}} &&\pi A_j\leq \bar{c}_j, j\in J, \label{RIOVLP_2-condition-1}\\
		&&\pi A_j= \bar{c}_j, j\in \bar{J}, \label{RIOVLP_2-condition-2}\\
		&&\bar{c}x^0=K. \label{RIOVLP_2-condition-3}
		\nqa
		then $\bar c^*$ is an optimal solution of the problem (\textit{RIOVLP}$_1$). 
	\end{thm}
	\begin{proof}
		Let $\mathcal{F}$ and $\mathcal{F}_1$ be the feasible regions of the problems (\textit{RIOVLP}$_1$) and (\textit{RIOVLP}$_{1}^1$), respectively. Notice that the constraints (\ref{RIOVLP_w-1})-(\ref{RIOVLP_w-2}) mean that $x^0$ is an optimal solution of the problem (\textit{LP}$_{\bar c}$), which holds if and only if its dual problem (\textit{DLP}$_{\bar{c}}$) has a feasible solution $\pi$ which satisfies the \textbf{complementary slackness conditions} $x^0_j(\bar{c}_j-\pi A_j)=0$ for any $j\in J\cup \bar J$. Suppose $\bar{c}\in \mathcal{F}$, then there exists ${\pi}$ satisfying the constraints  (\ref{RIOVLP_2-condition-1})-(\ref{RIOVLP_2-condition-2}).  Hence, we have $({\pi},\bar{c})\in \mathcal{F}_1$. On the other hand, suppose $({\pi}^*,\bar{c}^*)\in \mathcal{F}_1$ is an optimal solution of the problem (\textit{RIOVLP}$_{1}^1$), then $({\pi}^*,\bar{c}^*)$ satisfies the constraints (\ref{RIOVLP_2-condition-1})-(\ref{RIOVLP_2-condition-3}), which renders that $\bar{c}^*$ satisfies the constraint (\ref{RIOVLP_w-1})-(\ref{RIOVLP_w-2}). Hence, we have  $\bar{c}^*\in \mathcal{F}$. Furthermore, $\bar{c}^*$ is an optimal solution of the problem (\textit{RIOVLP}$_1$), since the two problems have the same optimal objective value. \hb
	\end{proof}
	
	Let $\bar{c}_j=c_j+\alpha_j-\beta_j$, where $\alpha_j,\beta_j$ are the increment and decrement of $c_j$,  respectively. We claim that at least one of $\alpha_j$ and $\beta_j$ is 0 based on the property of weighted $l_1$ norm for any $j\in J\cup \bar{J}$. 
	Thus, the problem (\textit{RIOVLP}$^1_{1}$) can be turned into the model below.
	\eqa
	\min &&\sum\limits_{j\in J\cup \bar  J}d_j(\alpha_j+\beta_j)\nn\\
	(\textit{RIOVLP}_{1}^2) \ {\hbox {s.t.}} &&\pi A_j\leq c_j+\alpha_j-\beta_j, j\in J, \label{URIOVLP_1-condition-1}\\
	&&\pi A_j= c_j+\alpha_j-\beta_j, j\in \bar{J}, \label{URIOVLP_1-condition-2}\\
	&&\sum\limits_{j\in \bar{J}}(c_j+\alpha_j-\beta_j)x^0_j=K, \label{URIOVLP_1-condition-3}\\
	&&\alpha_j \geq 0, j\in J\cup \bar{J},\label{URIOVLP_1-condition-4}\\
	&&\beta_j \geq 0, j\in J\cup \bar{J}. \label{URIOVLP_1-condition-5}
	\nqa
	Associate a dual variable $y_j$ with the constraints (\ref{URIOVLP_1-condition-1}) and (\ref{URIOVLP_1-condition-2}), and a dual variable $z$ with the constraint (\ref{URIOVLP_1-condition-3}). Then we can get its dual problem below.
	\eqa \label{RIOVLP_2-D}	
	\max && \sum_{j\in J\cup \bar{J}}c_j y_j+\Big{(}K-\sum_{j\in \bar{J}}c_jx^0_j\Big{)}z \nn\\
	(\textit{DRIOVLP}_1) \ {\hbox {s.t.}}  && Ay=0, \\
	&&  -y_j\leq d_j, j\in J, \\
	&&  -y_j+x^0_jz\leq d_j, j\in \bar{J},\\
	&&  y_j\leq d_j, j\in J,\\		 
	&&  y_j-x^0_jz\leq d_j, j\in \bar{J},\\
	&&  y_j\leq 0, j\in J.
	\nqa
	Delete the item $(K-\sum_{j\in \bar{J}}c_jx^0_j)z$ from the objective function of the problem (\textit{DRIOVLP}$_1$), we get a sub-problem $(D^z)$  below.
	\eqa \label{D-z}	
	\psi(z):= &\max& \sum_{j\in J\cup \bar{J}}c_j y_j  \nn\\
	(D^z) \ & {\hbox {s.t.}}  & Ay=0, \label{D-z-con-1}\\
	&&  -y_j\leq d_j, j\in J,  \label{D-z-con-2} \\
	&&  -y_j+x^0_jz\leq d_j, j\in \bar{J}, \label{D-z-con-3} \\
	&&  y_j\leq d_j, j\in J, \label{D-z-con-4} \\		 
	&&  y_j-x^0_jz\leq d_j, j\in \bar{J}, \label{D-z-con-5}\\
	&&  y_j\leq 0, j\in J. \label{D-z-con-6}
	\nqa
	Associate the dual variable $\bar\pi$ with the constraint (\ref{D-z-con-1}), $\bar\alpha_j$ with the constraints (\ref{D-z-con-2} and \ref{D-z-con-3}), $\bar\beta_j$ with the constraints (\ref{D-z-con-4}) and (\ref{D-z-con-5}). 
	If we treat $z$ as a constant, then the dual problem of $(D^z)$ can be stated as follows.
	\eqa
	\phi(z):=&\min &\sum_{j\in J\cup\bar{J}}d_j(\bar\alpha_j+\bar\beta_j)+z\sum_{j\in \bar{J}}(\bar\beta_j-\bar\alpha_j)x^0_j\nn\\
	(P^z) \ & {\hbox {s.t.}} &\bar\pi A_j-\bar\alpha_j+\bar\beta_j\leq c_j, j\in J,  \label{P-z-con-1}\\
	&& \bar\pi A_j-\bar\alpha_j+\bar\beta_j= c_j, j\in \bar{J}, \label{P-z-con-2}\\
	&& \bar\alpha_j \geq 0, j\in J\cup \bar{J}, \\ 
	&& \bar\beta_j \geq 0, j\in J\cup \bar{J}.
	\nqa
 Notice that the problem $(P^z)$ can be regarded as a Lagrange relaxation of (\textit{RIOVLP}$^2_{1}$) by introducing a Lagrange multiplier $z$ to the constraint (\ref{URIOVLP_1-condition-3}). 

Furthermore, let $y=x^0-\breve{y}$, then the problem $(D^z)$ can be turned into the following form.
		\eqa \label{D-bar-z}	
		\min && \sum_{j\in J\cup \bar{J}}c_j \breve{y}_j-\sum_{j\in  \bar{J}}c_j x^0_j  \nn\\
		(\breve D^z) \ {\hbox {s.t.}}  && A\breve{y}=b,  \\
		&&  0\leq \breve{y}_j\leq d_j, j\in J,    \\
		&&  (1-z)x^0_j-d_j\leq \breve{y}_j \leq (1-z)x^0_j+d_j, j\in \bar{J}.   \nqa
It is obvious that the problem $(\breve D^z)$  has the same objective function as (\textit{LP}) in (\ref{eq-LP}) by deleting the constant $-\sum_{j\in  \bar{J}}c_j x^0_j$.
		As for the feasible region, they have the same constraint condition 
		$A\breve{y}=b $. But they have different upper and lower bounds of the variables for a given $z$.  
		Hence, it provides us a way to solve $(\breve D^z)$ similar to the original problem (\textit{LP}). 
  
	
	\subsection{Properties of the function $\psi(z)$}
	
	Let $\mathcal{F}_{D^z}$ be the feasible region of the problem $(D^z)$ for a given $z\in \mathbb{R}$. 
	Next, we will discuss some properties of the feasible region $\mathcal{F}_{D^z}$ and objective  function $\psi(z)$ for a given $z\in \mathbb{R}$. 
	
	\begin{itemize}
		\item[$(P1)$] For a given $z\in \mathbb{R}$, if $\mathcal{F}_{D^z}\neq \emptyset$, then the problem $(D^z)$ has an optimal solution and $|\psi(z)|<+\infty$. It follows from $-d_j\leq y_j\leq 0$ for any $j\in J$ and $x^0_jz-d_j\leq y_j\leq x^0_jz+d_j$ for any $j\in \bar{J}$ by constraints (\ref{D-z-con-2})-(\ref{D-z-con-6}). 
		\item[$(P2)$] For any $z\in [\max_{j\in \bar{J}}\frac{-d_j}{x^0_j},\min_{j\in \bar{J}}\frac{d_j}{x^0_j}]$, we have $\bm{y}=\bm{0}\in \mathcal{F}_{D^z}\neq \emptyset$ and $|\psi(z)|<+\infty$.     
	\end{itemize}

	\begin{thm} \label{th-phi-property}
		Suppose $z_1<z_2$ and $|\psi(z_1)|,|\psi(z_2)|<+\infty$, then 
		the function $\psi(z)$ has the following properties for any $z\in [z_1,z_2]$.
		
		(1) $|\psi(z)|<+\infty$.
		
		(2) $\psi(z)$ is a continuous and piecewise function of $z$.
		
		
		(3) $\psi(z)$ is a concave function of $z$.
		
	\end{thm}
	\begin{proof}
We first claim that the assumption	$|\psi(z_1)|,|\psi(z_2)|<+\infty$ for $z_1<z_2$ must hold due to property $(P2)$.
		(1) Suppose $\bm{\tilde{y}}$ and $\bm{\bar{y}}$ are the optimal solutions of the problems $(D^{z_1})$ and $(D^{z_2})$ with objective values $\psi(z_1)$ and $\psi(z_2)$, respectively. Let $z=\lambda z_1+(1-\lambda)z_2$, where $0\leq \lambda \leq 1$. Next, we will prove $\lambda \bm{\tilde{y}}+(1-\lambda)\bm{\bar{y}}$ is a feasible solution of the problem $(D^z)$. It follows from $A\tilde{y}=0$ and $A\bar{y}=0$ that $A(\lambda \bm{\tilde{y}}+(1-\lambda)\bm{\bar{y}})=0$ holds. 
  It follows from $-\tilde{y}_j\leq d_j$ and $-\bar{y}_j\leq d_j$ that $-(\lambda \tilde{y}_j+(1-\lambda)\bar{y}_j)\leq d_j$ holds for any $j\in J$. Hence, the constraint (\ref{D-z-con-2}) holds. In a similar way, the constraints (\ref{D-z-con-4}) and (\ref{D-z-con-6}) hold. It follows from $-\tilde{y}_j\leq d_j-x^0_jz_1$ and $-\bar{y}_j\leq d_j-x^0_jz_2$ that $-(\lambda \tilde{y}_j+(1-\lambda)\bar{y}_j)\leq \lambda (d_j-x^0_jz_1)+(1-\lambda)(d_j-x^0_jz_2)=d_j-x^0_j(\lambda z_1+(1-\lambda) z_2)=d_j-x^0_jz$ for any $j\in \bar{J}$. Hence, the constraint (\ref{D-z-con-3}) holds. In a similar way, the constraint (\ref{D-z-con-5}) holds. Therefore, we have $\mathcal{F}_{D^z}\neq \emptyset$ and $|\psi(z)|<+\infty$ for any $z\in [z_1,z_2]$ by property $(P1)$.
		
		(2)  We first claim that $\psi(z)=\phi(z)$ for any $z\in [z_1,z_2]$ based on the strong duality of linear programming.  
		Furthermore, the problem $(P^z)$ has an optimal solution $(\bar\pi^z,\bar\alpha^z,\bar\beta^z)$ since its dual problem $(D^z)$ is feasible and bounded for any $z\in [z_1,z_2]$ by (1).  
		Therefore, there exists an optimal solution $(\bar\pi^z,\bar\alpha^z,\bar\beta^z)$ satisfying $|\bar\pi^z_i|,|\bar\alpha^z_j|,|\bar\beta^z_j|\leq M$ for any $z\in [z_1,z_2]$, where $M=n!c_{max}$ with $c_{max}=\max_{j=1,\cdots,n} |c_j|$ is a large positive number by Lemma 2.1 in \cite{Papadimitriou_book_1982}. Hence, if we add a set of constraints: $|\bar\pi_i|\leq M$, $|\bar\alpha_j|\leq M$ and $|\bar\beta_j|\leq M$ to the problem $(P^z)$, then the optimal solution of problem $(P^z)$ will be unchanged. It renders that the feasible region of $(P^z)$  becomes a bounded closed convex set. Hence, $\phi(z)$, as well as $\psi(z)$, is a continuous and piecewise function of $z$ for $z\in [z_1,z_2]$.

		
		(3) Suppose $z_l<z_r$ and $z_l,z_r\in [z_1,z_2]$. 
		Let $(\bar\pi^*,\bar\alpha^*,\bar\beta^*)$ be an optimal solution of the problem $(P^z)$, where $z=kz_l+(1-k)z_r$ and $0\leq k \leq 1$. Then we have   
		\eqa 
		&&\phi(z)=\phi(kz_l+(1-k)z_r)\nn\\
		&=&\sum_{j\in J\cup\bar{J}}d_j(\bar\alpha^*_j+\bar\beta^*_j)+(kz_l+(1-k)z_r)\sum_{j\in \bar{J}}(\bar\beta^*_j-\bar\alpha^*_j)x^0_j\nn\\
		&=&k\Bigg(\sum_{j\in J\cup\bar{J}}d_j(\bar\alpha^*_j+\bar\beta^*_j)+z_l\sum_{j\in \bar{J}}(\bar\beta^*_j-\bar\alpha^*_j)x^0_j\Bigg)\nn\\
		&+&(1-k)\Bigg(\sum_{j\in J\cup\bar{J}}d_j(\bar\alpha^*_j+\bar\beta^*_j)+z_r\sum_{j\in \bar{J}}(\bar\beta^*_j-\bar\alpha^*_j)x^0_j\Bigg)\nn\\
		&\geq&k\phi(z_l)+(1-k)\phi(z_r).
		\nn \nqa
		Hence, $\phi(z)$,  as well as $\psi(z)$,  is a concave function of $z$.  \hb
	\end{proof}
	
	Now we claim that the feasible region $\mathcal{F}_{D^z}$ may be empty for a given $z\in \mathbb{R}$. Next we present an example to illustrate it. Let $A=[E_m,E_m]$ and $z_0=\frac{2d_{max}}{x^0_{min}}$, where $E_m$ is an identity matrix, $x^0_{min}=\min_{j\in \bar{J}}x^0_j$ and $d_{max}=\max_{j\in J\cup \bar{J}}d_j$. 
	Then $x^0_jz-d_j\geq x^0_jz-d_{max}> x^0_j\frac{2d_{max}}{x^0_{min}}-d_{max}\geq 2d_{max}-d_{max}=d_{max}$ for any $j\in \bar{J}$ when $z>z_0$. Therefore, we have $d_{max}<x^0_jz-d_j\leq y_j\leq x^0_jz+d_j$ for any $j\in \bar{J}$. Consider the $i$-th constraint in (\ref{D-z-con-1}), we have $y_i+y_{i+m}=0$. If $i\in \bar{J}$ or $i+m\in \bar{J}$, then this constraint can not hold since $-d_j\leq y_j\leq 0$ for any $j\in J$ and $y_j>d_{max}$ for any $j\in \bar{J}$. Unfortunately, this case must exist since $\bar{J}\neq \emptyset$. Hence, $\mathcal{F}_{D^z}$ may be an empty set for a given $z>0$. In a similar way, $\mathcal{F}_{D^z}$ may be an empty set for a given $z<0$. In these cases, $\psi(z)$ is undefined.
	
	\begin{thm} \label{th-Continuity-of-Empty-Sets}
		(1) If $z_1<0$ and $\mathcal{F}_{D^{z_1}}=\emptyset$, then $\mathcal{F}_{D^{z}}=\emptyset$ for any $z\leq z_1$.
		
		(2) If $z_2>0$ and $\mathcal{F}_{D^{z_2}}=\emptyset$, then  $\mathcal{F}_{D^{z}}=\emptyset$ for any $z\geq z_2$.		
	\end{thm}
	\begin{proof}
		(1) Suppose there exists $z_0<z_1$ satisfying $\mathcal{F}_{D^{z_0}}\neq \emptyset$. Then we have $\mathcal{F}_{D^{z}}\neq \emptyset$ for any $z\in [z_0,0]$ by Theorem \ref{th-phi-property}, which contradicts that $\mathcal{F}_{D^{z_1}}=\emptyset$. Hence, the conclusion holds. In a similar way, (2) holds. \hb
	\end{proof}

 For simplicity of discussion, we turn the problem $(D^z)$ into a standard \textit{LP} problem by the following steps. 
	(i) It follows from the constraint (\ref{D-z-con-6}) and $d_j>0$ that (\ref{D-z-con-4}) can be omitted. (ii) Add relaxation variables $\xi_j$ to the constraints (\ref{D-z-con-2}) and (\ref{D-z-con-3}). (iii) Add relaxation variables $\eta_j$ to the constraints (\ref{D-z-con-5}) and (\ref{D-z-con-6}). (iv) Replace non-constrained variables $y_j$ by two non-negative variables $\dot{y}_j$ and $\ddot{y}_j$. Then we get its maximization standard \textit{LP} problem below.
	\eqa \label{D-z-s}	
	\psi(z):=&\max & \sum_{j\in J\cup \bar{J}}c_j(\dot{y}_j-\ddot{y}_j)  \nn\\
	(D^z_s) \ & {\hbox {s.t.}}  & A\dot{y}-A\ddot{y}=0,  \\
	&&  -\dot{y}_j+\ddot{y}_j+\xi_j= d_j, j\in J,   \label{eq-Dzs-con2} \\ 
	&&  -\dot{y}_j+\ddot{y}_j+\xi_j= d_j-x^0_jz, j\in \bar{J},   \label{eq-Dzs-con3}  \\	 
	&&  \dot{y}_j-\ddot{y}_j+\eta_j= 0, j\in J,  \\
	&&  \dot{y}_j-\ddot{y}_j+\eta_j= d_j+x^0_jz, j\in \bar{J},  \label{eq-Dzs-con5}  \\
	&& \dot{y},\ddot{y},\xi,\eta \geq 0. 
	\nqa
	Let $\bar{A}$, $\bar{b}^z$, $\bar{c}$ and $\bar{X}$ be the coefficient matrix, right-hand vector, cost coefficient vector and all variables of the problem $(D^z_s)$, respectively. For the convenience of further discussion, we assume  $J=\{1,2,\cdots,k\}$ and $\bar{J}=\{k+1,\cdots,n\}$, where $0\leq k \leq n-1$.  Then we have 
	\eq \bar{A}=\begin{bmatrix}
		A & -A & \bm{0}  &  \bm{0} \\
		-E_n & E_n & E_n & \bm{0} \\
		E_n & -E_n & \bm{0} & E_n \\
	\end{bmatrix}_{\bar{m} \times  \bar{n}}, \label{eq-barA}\nq
	$\bar{c}=(c,-c,\bm{0},\bm{0})$, $\bar{b}^z=(\bm{0},d_J,d_{\bar{J}}-zx^{0}_{\bar{J}},\bm{0},d_{\bar{J}}+zx^{0}_{\bar{J}})^T$, $\bar{X}^T=(\dot{y},\ddot{y},\xi,\eta)$, where $\bar{m}=2n+m, \bar{n}=4n$ and $x^{0}_{\bar{J}}$ is a row vector here. The corresponding standard \textit{LP} model can be formulated below.
 \eqa \label{D-z-s1}	
	\psi(z):=&\max & \bar c \bar X  \nn\\
	(D^z_s) \ & {\hbox {s.t.}}  & \bar A\bar X=\bar b^z, \nn \\
&&\bar X \geq 0. \nn
	\nqa
	
	Next we need to clarify the maximum value $z_l$ and minimum  value  $z_r$ satisfying the property (1) and (2) in Theorem \ref{th-Continuity-of-Empty-Sets}, respectively. For convenience, we call them left and right break points, respectively.
	
	\begin{defn}
		A value $z_l < 0$ is called \textbf{the left break point} of $\psi(z)$, if $\mathcal{F}_{D^{z}}=\emptyset$ for any $z< z_l$ and $\mathcal{F}_{D^{z}}\neq \emptyset$ for any $z\in [z_l,0]$. A value $z_r > 0$ is called \textbf{the right break point} of $\psi(z)$, if $\mathcal{F}_{D^{z}}=\emptyset$ for any $z>z_r$ and $\mathcal{F}_{D^{z}}\neq \emptyset$ for any $z\in [0,z_r]$.
	\end{defn}
	
	\begin{thm}\label{th-breakpoint}
		(1)  If there is a value $z_1<0$ satisfying $\mathcal{F}_{D^{z_1}}=\emptyset$, then there exists a left break point $z_l$.  
		(2)  If there is a value $z_2>0$ satisfying  $\mathcal{F}_{D^{z_2}}=\emptyset$, then there exists a right break point  $z_r$. (3) For any $z \geq z_l$ and $z \leq z_r$, $|\psi(z)|<+\infty$. 
	\end{thm}
	\begin{proof}
		(1) Note that $\mathcal{F}_{D^z}\neq \emptyset$ for any $z\in [\max_{j\in \bar{J}}\frac{-d_j}{x^0_j},\min_{j\in \bar{J}}\frac{d_j}{x^0_j}]$ by property (P2). Hence, we can assume that there exists $z_0< 0$ satisfying $\mathcal{F}_{D^{z_0}}\neq \emptyset$. 
		Let $B$ be the optimal basis corresponding to the problem $(D^{z_0}_s)$.  
		Consider the $i$-th basic variable:
		$\bar{X}_i^{z_0}=B^{-1}_i\bar{b}^{z_0}=k_iz_0+g_i$, where $B^{-1}_i$ is the $i$-th row of the inverse matrix $B^{-1}$ and $k_i,g_i$ are constants. Therefore, when the optimal basis $B$ changes to a non-optimal basis as $z_0$ decreases, there must be a non-negative basic variable that changes to a negative value. Suppose $z_l$ is such a value. Then we have $\mathcal{F}_{D^{z_l}}\neq \emptyset$. Furthermore, it follows from property (1) in Theorem \ref{th-phi-property} that $\mathcal{F}_{D^{z}}\neq \emptyset$ for any $z\in [z_l,0]$. Hence, $z_l$ is the left break point.  The property (2) holds similarly. (3) This property holds due to (1) and (2). \hb

		\end{proof}

\begin{figure}[htbp] 
	\centering	
 \includegraphics[totalheight=6.5cm]{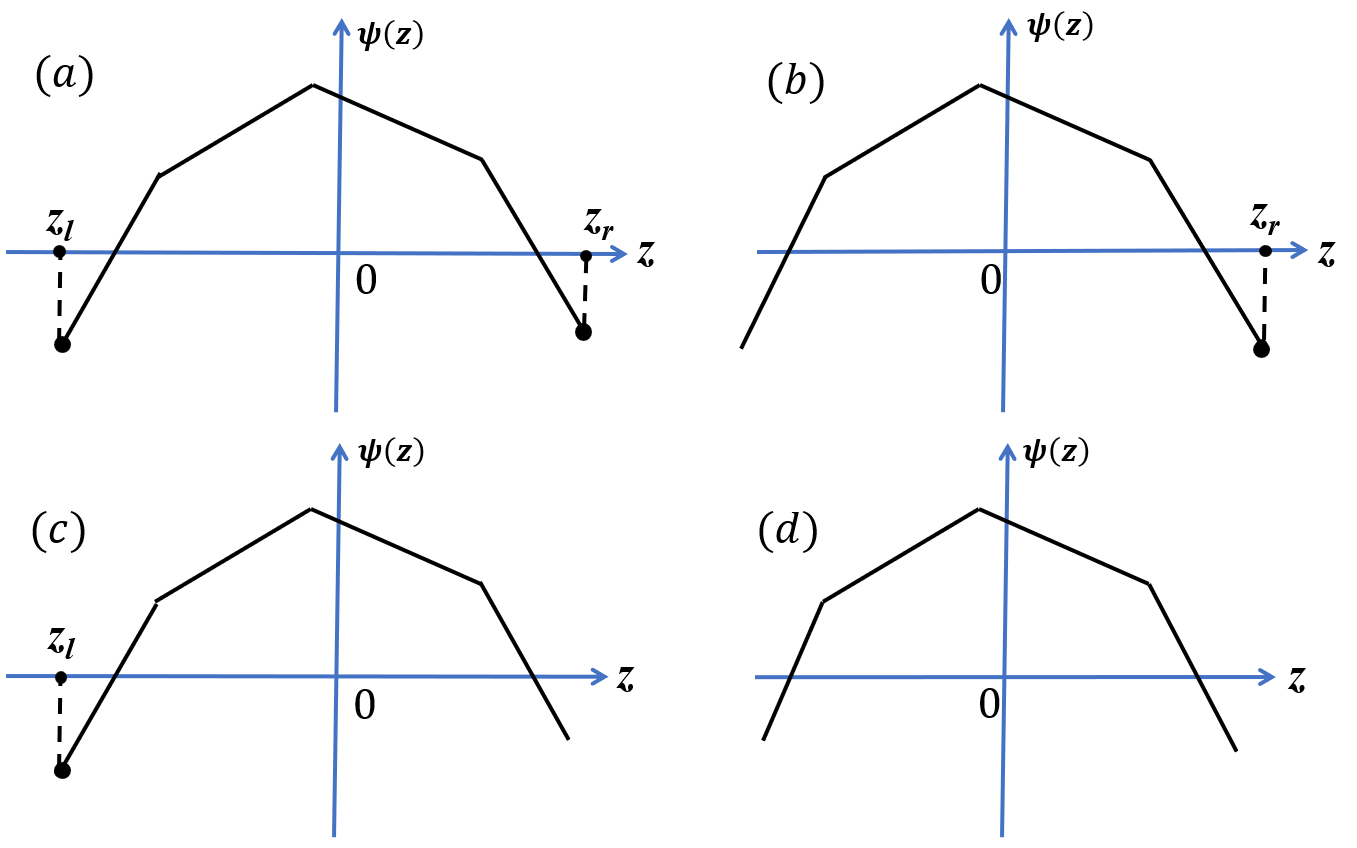}
\caption{
	Four cases of left and right break points in the function $\psi(z)$.\\  $(a)$ $z_l,z_r$ exist; $(b)$ $z_r$ exists; $(c)$ $z_l$ exists; $(d)$ no $z_l,z_r$.
 }	\label{figure-phi}	
		\end{figure}

		
		Based on the previous discussions, we claim that there are four cases of the left and right break points in the function $\psi(z)$ as shown in Figure \ref{figure-phi}, where $z_l$ and $z_r$ represent the left and right break points, respectively.
Notice that $\psi(z)$ is a continuous piecewise linear function, we introduce the following definition.

		\begin{defn}
			Let $(z_i, \psi(z_i))$ be an intersection point of two adjacent lines. We call $z_i$ \textbf{a turning coordinate}.
		\end{defn}
  
  \section{Solve the problem (\textit{RIOVLP}$_1$) when A is unimodular} \label{section_RIOVLP_unimodular}
  
  In this section, we solve the problem (\textit{RIOVLP}$_1$) when A is a unimodular coefficient matrix. Firstly we analyze some properties of turning coordinates. Secondly, we  present an important theorem to determine the critical value $z^*$ in an optimal solution ($y^*,z^*$) of the problem (\textit{DRIOVLP}$_1$). Thirdly, we calculate the slope $k_z$ of a piece of segment in $\psi(z)$ for a given $z$ which is not a turning coordinate. Then we present two algorithms to calculate the left and right break points $z_l$ and $z_r$. Finally, we propose an algorithm to solve the problem (\textit{RIOVLP}$_1$).
		\begin{asm} \label{assumption-A}
			The coefficient matrix $A$ is unimodular.
		\end{asm}
		\begin{asm} \label{assumption-x^0}
			The given feasible solution $x^0$ and the weight vector $d$ are integral.
		\end{asm}

  \subsection{Properties of turning coordinates}
		To describe the form of left and right break points, as well as turning points, we present the following lemma.
		\begin{lem} \label{lemma-TC-form}
	Suppose Assumptions \ref{assumption-A} and \ref{assumption-x^0} hold.		If $z_i$ is a turning coordinate, then $z_i$ is in the form of $\frac{\varsigma}{\sigma}$, where $\varsigma,\sigma$ are integers and $|\sigma|\leq 2|\bar{J}|x^0_{max}$, $|\varsigma|\leq (n+|\bar{J}|)d_{max}$, $x^0_{max}=\max_{j\in \bar{J}}x^0_j$. 
		\end{lem}
		\begin{proof} 
			It follows from Assumption \ref{assumption-A} that the coefficient matrix $\bar{A}$ defined as in (\ref{eq-barA}) of the problem $(D^z_s)$ is also a unimodular matrix. For a given $z\in \mathbb{R}$, if $\mathcal{F}_{D^z}\neq \emptyset$, then there exists an optimal basis matrix $B$ of the problem $(D^z_s)$. Hence, the elements $b^-_{ij}\in \{0,1,-1\}$ for the inverse matrix $B^{-1}=(b^-_{ij})$. Specifically, let $b^1_{ij}, b^2_{ij}, b^3_{ij}$ be the elements of $B_i^{-1}$ corresponding to $j$ in (\ref{eq-Dzs-con2}), (\ref{eq-Dzs-con3}) and (\ref{eq-Dzs-con5}), respectively.
			Let $\bar{X}_i^z$ be the $i$-th basic variable corresponding to $B$. Then we have
			\eqa
			\bar{X}_i^z&=& B^{-1}_i\bar{b}^z \nn\\
			&=&\sum_{j\in J}b^1_{ij}d_j+\sum_{j\in \bar{J}}b^2_{ij}(d_j-x^0_jz)+\sum_{j\in \bar{J}}b^3_{ij}(d_j+x^0_jz) \nn\\
			&=& \sum_{j\in J}b^1_{ij}d_j+\sum_{j\in\bar{ J}}(b^2_{ij}d_j+b^3_{ij}d_j)+z\sum_{j\in \bar{J}}x^0_j(b^3_{ij}-b^2_{ij}).\nn
			\nqa
			
			With the variation of $z$, the optimal basis matrix $B$ of the problem $(D^z_s)$ will change at a value in the form of $z=\frac{\sum_{j\in J}b^1_{ij}d_j+\sum_{j\in\bar{J}}(b^2_{ij}d_j+b^3_{ij}d_j)}{\sum_{j\in \bar{J}}x^0_j(b^2_{ij}-b^3_{ij})}$ as in the proof of Theorem \ref{th-breakpoint}.
			It is easy to know $|\sum_{j\in J}b^1_{ij}d_j+\sum_{j\in\bar{J}}(b^2_{ij}d_j+b^3_{ij}d_j)|\leq (n+|\bar{J}|)d_{max}$ is an integer. Furthermore, $|\sum_{j\in \bar{J}}x^0_j(b^2_{ij}-b^3_{ij})|\leq 2|\bar{J}|x^0_{max}$ is also an integer. As a conclusion, the lemma holds.	
			\hb
		\end{proof}
		
		\begin{cor} \label{turning-coordinate-interval}
			Suppose $z_i$, $z_j$ are two turning coordinates, then $\displaystyle |z_i-z_j|\geq \frac{1}{(2|\bar{J}|x^0_{max})^2}$.
		\end{cor}
		\begin{proof}
			Suppose $z_i=\frac{\varsigma_i}{\sigma_i}$ and $z_j=\frac{\varsigma_j}{\sigma_j}$, then 
			$\varsigma_i,\varsigma_j,\sigma_i,\sigma_j\in \mathbb{Z}$ and $|\sigma_i|,|\sigma_j|\leq 2|\bar{J}|x^0_{max}$. Hence,  
			\eq \left|\frac{\varsigma_i}{\sigma_i}-\frac{\varsigma_j}{\sigma_j}\right| = 
			\left| \frac{\varsigma_i\sigma_j-\varsigma_j\sigma_i}{\sigma_i\sigma_j}\right|  \geq \frac{\left |\varsigma_i\sigma_j-\varsigma_j\sigma_i \right|}{(2|\bar{J}|x^0_{max})^2}\geq \frac{1}{(2|\bar{J}|x^0_{max})^2}. \nn
			\nq
This completes the proof.			\hb 
		\end{proof}
		
		Let $TC=\{\frac{\varsigma}{\sigma}|\varsigma,\sigma\in \mathbb{Z}, |\sigma|\leq 2|\bar{J}|x^0_{max}, |\varsigma|\leq (n+|\bar{J}|)d_{max}\}$, then a turning coordinate must in $TC$, but the opposite is not true. 
		
		It follows from Lemma \ref{lemma-TC-form} that the optimal basis of the problem $(D^{z_0}_s)$ keeps unchanged in its $\epsilon$-neighborhood for a given value $z_0\notin TC$ and $\epsilon>0$. It renders that the optimal solution $(\bar\pi^*,\bar\alpha^*,\bar\beta^*)$ of the problem $(P^{z_0})$ keeps unchanged for $z\in [z_0-\epsilon,z_0+\epsilon]$. 
		Therefore, the slope of $\psi(z)$ is $k_z=\sum_{j\in \bar{J}}(\bar\beta^*_j-\bar\alpha^*_j)x^0_j$ for $z\in [z_0-\epsilon,z_0+\epsilon]$.
		\begin{cor} \label{cor-continuity of basis}
			For a given value $z_0\notin TC$, if $\mathcal{F}_{D^{z_0}}\neq \emptyset$ and $B$ is an optimal basis of the problem $(D^{z_0}_s)$, then there must exist $\epsilon>0$ satisfying that $B$ is also an optimal basis of the problem $(D^z_s)$ for any $z\in [z_0-\epsilon,z_0+\epsilon]$. 
		\end{cor}

		\subsection{Solve the dual inverse problem (\textit{DRIOVLP}$_1$)}
		For the convenience of following discussion, we give the definition below.
		
		\begin{defn} \label{critical-value}
			Suppose $(y^*,z^*)$ is an optimal solution of the problem (\textit{DRIO} \textit{VLP}$_1$), then we call $z^*$ \textbf{a critical value}.
		\end{defn}
		
		If we can determine  \textbf{a critical value} $z^*$, then we can solve the sub-problem $(D^{z^*})$ and obtain its optimal solution $y^*$, which renders an optimal solution $(y^*,z^*)$ of the problem (\textit{DRIOVLP}$_1$). Hence, the main issue now is how to determine $z^*$. To do so, we need to determine two critical  consecutive segments in the piecewise linear function $\psi(z)$, whose slopes $k_1, k_2$ satisfy $k_2 < c_{\bar J}x^0_{\bar J}-K\leq k_1$.
		Notice that there is only one slope of a linear function for the left and right break points and there are cases that the left and right break points do not exist as shown in Figure \ref{figure-phi}.  We first introduce the minimum and maximum turning coordinates  $\underline z$ and $\bar{z}$ to substitute the the left and right break points. Then  define the right slope of $\psi(z)$ (if exists) as $k^+_z$ and the left slope of $\psi(z)$ (if exists) as $k^-_z$.
		
\begin{figure}[htbp]  
			\centering
			\includegraphics[totalheight=6.5cm]{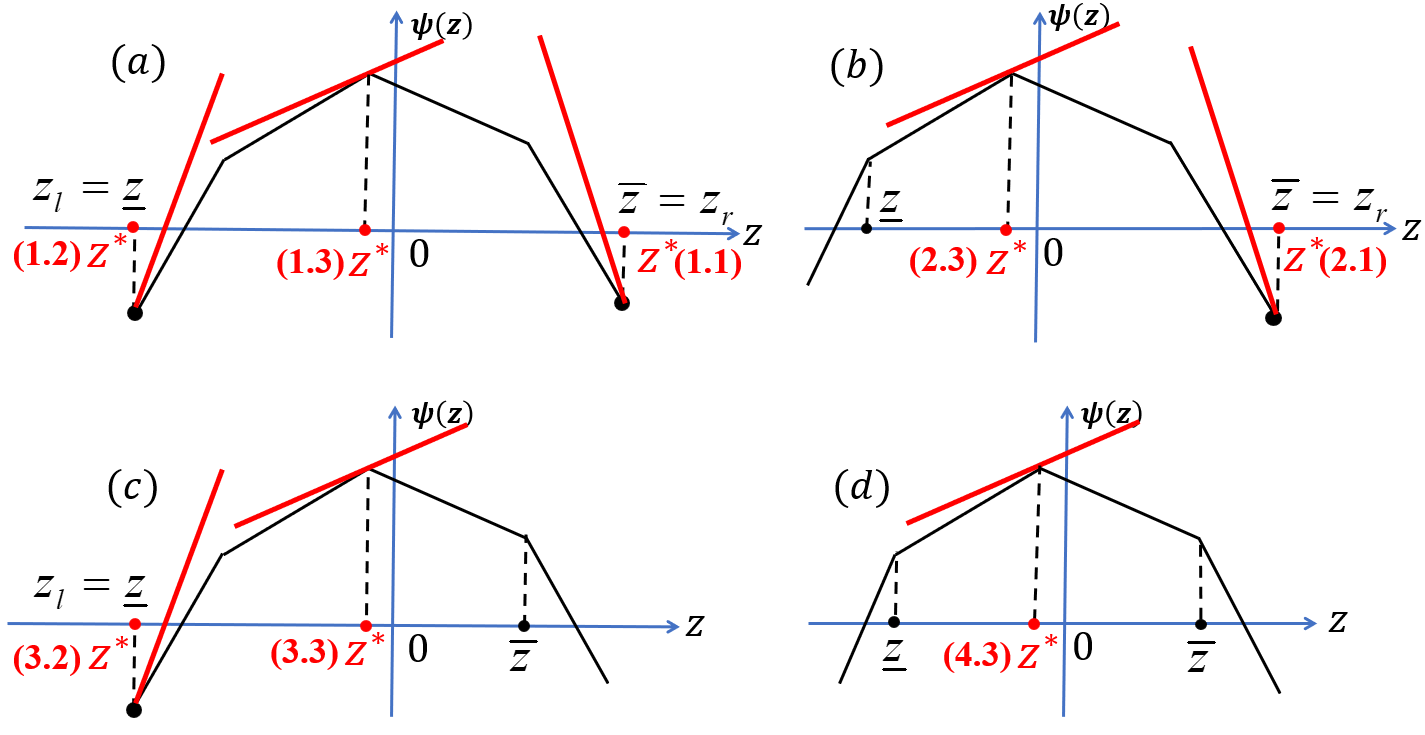}
			\caption{The values $\underline z$, $\bar{z}$ and $z^*$ in each case of Theorem \ref{th-z*}.}	\label{figure-phi-optimality}	
		\end{figure}

		Next, we propose an optimality condition to determine $z^*$ at each of the four cases. Figure \ref{figure-phi-optimality} shows the values $\underline z$, $\bar{z}$ and $z^*$ in each case, where the red line represents the line with the slope $\delta=c_{\bar J}x^0_{\bar J}-K$ in each subcases. 
		
		\begin{thm}\label{th-z*}
			Suppose $(\hat z,\psi(\hat z))$ is an intersection point of two consecutive segments in the piecewise linear function $\psi(z)$.
			
			{\rm (1)} Suppose $\mathcal{F}_{D^z}=\emptyset$ for any $z\in (-\infty,\underline{z})\cup (\bar{z},+\infty)$. Then 
			
			{\rm (1.1)} If $k^-_{{\bar{z}}}\geq \delta$, then $z^*={\bar{z}}$ is \textbf{the critical value}. 
			
			{\rm (1.2)} If  $k^+_{{\underline{z}}}<\delta$, then $z^*={\underline{z}}$ is \textbf{the critical value}.  
			
			{\rm (1.3)} If $k_{\hat z}^+ < \delta\leq k_{\hat z}^-$, then $z^*=\hat z$ is \textbf{the critical value}. 
			
			{\rm (2)} Suppose $\mathcal{F}_{D^z}\neq\emptyset$ for any $z\in (-\infty,\underline{z})$  and $\mathcal{F}_{D^z}=\emptyset$ for any $z\in (\bar{z},+\infty)$. Then 
			
			{\rm (2.1)} If $k^-_{{\bar{z}}}\geq \delta$, then $z^*={\bar{z}}$ is \textbf{the critical value}. 
			
			{\rm (2.2)} If  $k^-_{{\underline{z}}}<\delta$, then the problem (\textit{RIOVLP}$_1$) is infeasible. 
			
			{\rm (2.3)} If $k_{\hat z}^+ < \delta\leq k_{\hat z}^-$, then $z^*=\hat z$ is \textbf{the critical value}. 
			
			{\rm (3)} Suppose $\mathcal{F}_{D^z}=\emptyset$ for any $z\in (-\infty,\underline{z})$ and $\mathcal{F}_{D^z}\neq\emptyset$ for any $z\in (\bar{z},+\infty)$. Then 
			
			{\rm (3.1)} If $k^+_{{\bar{z}}}>\delta$, then the problem (\textit{RIOVLP}$_1$) is infeasible. 
			
			{\rm(3.2)} If  $k^+_{{\underline{z}}}\leq \delta$, then $z^*={\underline{z}}$ is \textbf{the critical value}.  
			
			{\rm(3.3)} If $k_{\hat z}^+ \leq \delta< k_{\hat z}^-$, then $z^*=\hat z$ is \textbf{the critical value}. 
			
			{\rm (4)} Suppose $\mathcal{F}_{D^z}\neq\emptyset$ for any $z\in (-\infty,\underline{z})\cup (\bar{z},+\infty)$. Then 
			
			{\rm(4.1)} If $k^+_{{\bar{z}}}>\delta$, then the problem (\textit{RIOVLP}$_1$) is infeasible. 
			
			
			{\rm (4.2)} If  $k^-_{{\underline{z}}}<\delta$, then the problem (\textit{RIOVLP}$_1$) is infeasible. 
			
			{\rm (4.3)} If $k_{\hat z}^+ \leq \delta\leq k_{\hat z}^-$, then $z^*=\hat z$ is \textbf{the critical value}. 
		\end{thm}
		
		\textbf{Proof.}
		(1) (1.1) Suppose $(\tilde{y},\tilde{z})$ and $\bar{y}$ are optimal solutions of the problems (\textit{DRIOVLP}$_1$) and $(D^{\bar{z}})$, respectively. It is easy to know that $\tilde{z} \leq \bar{z}$. 
		If $\tilde{z}<\bar{z}$,  then 
		\eqa &&\sum_{j\in J\cup \bar{J}}c_j\tilde{y}_j+(K-\sum_{j\in \bar{J}}c_jx^0_j)\tilde{z}= c\tilde{y} -\delta\tilde{z}  \nn\\
		&\leq& \psi(\tilde{z})-\delta\tilde{z}\nn\\
		& \leq & k^-_{\bar{z}}(\tilde{z}-\bar{z})+\psi(\bar{z})-\delta\tilde{z},\  (\psi(z) \ is \ concave\ by\ Theorem\ \ref{th-phi-property}.)\nn\\
		&= & \psi(\bar{z})-\delta\bar{z}+(\tilde{z}-\bar{z})(-\delta+k^-_{{\bar{z}}}) \nn\\ 
		&\leq & \psi({\bar{z}})- \delta{\bar{z}}\nn\\
		& = &\ c\bar{y} -\delta\bar{z}. \nn
		\nqa
		Hence, $(\bar{y},\bar{z})$ is also an optimal solution of the problem (\textit{DRIOVLP}$_1$) if the equalities always hold in the above three inequalities. Otherwise, it contradicts the optimality of $(\tilde{y},\tilde{z})$.
		
		
		(1.2) Suppose $(\tilde{y},\tilde{z})$ and $\underline{y}$ are optimal solutions of (\textit{DRIOVLP}$_1$) and $(D^{\underline{z}})$, respectively. It is easy to know that $\tilde{z} \geq \underline{z}$. If $\tilde{z}>\underline{z}$, then 
		\eqa &&c\tilde{y} -\delta\tilde{z} \nn\\
		&\leq& \psi(\tilde{z})-\delta\tilde{z}\nn\\
		& \leq & k^+_{\underline{z}}(\tilde{z}-\underline{z})+\psi(\underline{z})-\delta\tilde{z},\  (\psi(z) \ is \ concave\ by\ Theorem\ \ref{th-phi-property}.)\nn\\
		&= & \psi(\underline{z})-\delta\underline{z}+(\tilde{z}-\underline{z})(-\delta+k^+_{{\underline{z}}}) \nn\\ 
		&< & \psi({\underline{z}})-\delta{\underline{z}}\nn\\
		& = &c\underline{y} -\delta{\underline{z}}. \nn
		\nqa
		which contradicts that $(\tilde{y},\tilde{z})$ is an optimal solution of (\textit{DRIOVLP}$_1$).
		
		(1.3) Suppose $(\tilde{y},\tilde{z})$ and $\hat y$ are optimal solutions of  (\textit{DRIOVLP}$_1$) and $(D^{\hat z})$, respectively. We will discuss it in two situations. 
		(i) If $\tilde{z}>\hat z$, then 
		\eqa &&
  c\tilde{y} -\delta\tilde{z}\nn\\
		&\leq& \psi(\tilde{z})-\delta\tilde{z}\nn\\
		& \leq & k^+_{\hat z}(\tilde{z}-\hat z)+\psi(\hat z)-\delta\tilde{z},\  (\psi(z) \ is \ concave\ by\ Theorem\ \ref{th-phi-property}.)\nn\\
		&= & \psi(\hat z)-\delta\hat z+(\tilde{z}-\hat z)(-\delta+k_{\hat z}^+) \nn\\ 
		&< & \psi({\hat z})-\delta{\hat z}\nn\\
		& = &c\hat{y} -\delta{\hat z}. \nn
		\nqa
		which contradicts that $(\tilde{y},\tilde{z})$ is an optimal solution of the problem (\textit{DRIOVLP}$_1$).
		
		(ii) If  $\tilde{z}<\hat z$, then 
		\eqa &&c\tilde{y} -\delta\tilde{z} \nn\\
		&\leq& \psi(\tilde{z})-\delta\tilde{z}\nn\\
		& \leq & k^-_{\hat z}(\tilde{z}-\hat z)+\psi(\hat z)-\delta\tilde{z},\  (\psi(z) \ is \ concave\ by\ Theorem\ \ref{th-phi-property}.)\nn\\
		&= & \psi(\hat z)-\delta\hat z+(\tilde{z}-\hat z)(-\delta+k_{\hat z}^-) \nn\\ 
		&\leq & \psi({\hat z})-\delta{\hat z}\nn\\
		& = &c\hat{y} -\delta{\hat z}. \nn
		\nqa  
		Then  $(\hat y,\hat z)$ is also an optimal solution of problem (\textit{DRIOVLP}$_1$) if the equalities always hold in the above three inequalities. Otherwise, it contradicts the optimality of $(\tilde{y},\tilde{z})$.
		
		Notice that one and only one case holds due to the concavity of $\psi(z)$.
		
		In a similar way, we can show that (2),(3) and (4) hold. \hb

		\subsection{Calculate the slope $k_z$ for a given $z\notin TC$ } \label{sub_section_cal_k_z}
		In this subsection, we focus on addressing the following two questions. For a given $z\in \mathbb{R}$, (Q1) if $\mathcal{F}_{D^{z}}\neq\emptyset$, then how to determine whether $z$ is a turning coordinate of an intersection point $(z,\psi(z))$ of two consecutive segments in $\psi(z)$? (Q2) if $z\notin TC$ and $\mathcal{F}_{D^{z}}\neq\emptyset$, then how to calculate the slope $k_z$ and obtain the  left and right  turning coordinates closest to $z$?  To solve these problems, we  make the assumption below.

		\begin{asm} \label{asm-solve-D-z}
			There is an algorithm $\mathcal{A}_z$ not only to determine whether $\mathcal{F}_{D^{z}}$ is an empty set for a given $z\in \mathbb{R}$ but also to solve the problem $(D^z)$ or $(D^z_s)$ or $(\breve{D}^z)$ and calculate $\psi(z)$ when $\mathcal{F}_{D^{z}}\neq\emptyset$.
		\end{asm}
		
		Let $\Delta=\frac{1}{(2|\bar{J}|x^0_{max})^2}$. 
		To answer the question (Q1), the following lemma is clearly valid.
		\begin{lem}
			Suppose Assumption \ref{asm-solve-D-z} holds and $\mathcal{F}_{D^{z}}\neq\emptyset$. Let $z_1=z-\frac{\Delta}{2}$, $z_2=z-\frac{\Delta}{4}$, $z_3=z+\frac{\Delta}{4}$,$z_4=z+\frac{\Delta}{2}$, $k_{z,z_1}=\frac{\psi(z)-\psi(z_1)}{z-z_1}$,$k_{z,z_2}=\frac{\psi(z)-\psi(z_2)}{z-z_2}$,$k_{z,z_3}=\frac{\psi(z)-\psi(z_3)}{z-z_3}$ and $k_{z,z_4}=\frac{\psi(z)-\psi(z_4)}{z-z_4}$. If $k_{z,z_1}=k_{z,z_2}$, $k_{z,z_3}=k_{z,z_4}$ and $k_{z,z_2}\neq k_{z,z_3}$, then $z$ is a turning coordinate of an intersection point $(z,\psi(z))$ of two consecutive segments in the function $\psi(z)$.
		\end{lem}

		
		Next, we concentrate on (Q2) calculating the slope $k_z$ for a given $z\notin TC$ and $\mathcal{F}_{D^{z}}\neq\emptyset$. 
		It seems straightforward to calculate $k_z=\frac{\psi(z+\epsilon)-\psi(z)}{\epsilon}$, where $\epsilon$  is a very small positive real number.
		However, it is a dilemma to determine $\epsilon$ such that $(z,z+\epsilon)$ does not contain a turning coordinate. 
		Hence, we calculate $k_z$ by duality theory of \textit{LP}. 
		We can obtain  an optimal solution $(\bar\pi^*,\bar\alpha^*,\bar\beta^*)$ of the problem $(P^z)$ based on an optimal solution $y^*$ of the problem $(D^z)$. 
		Therefore, we have $k_z=\sum_{j\in \bar{J}}(\bar\beta^*_j-\bar\alpha^*_j)x^0_j$ by Corollary \ref{cor-continuity of basis}. 
		
		If $\mathcal{F}_{D^{z_0}}\neq\emptyset$ and $z_0\notin TC$, then we will discuss how to calculate the left and right turning coordinates closest to $z_0$ in detail. First, we choose $z_1\in \mathbb{R}$ satisfying the conditions (C1), (C2) or (C1), (C3) below.
  
		(C1) $z_1 \neq z_0$ and $z_1\notin TC$. 
		(C2) If $z_1<z_0$, then $[z_1,z_0]\cap TC=\emptyset$. 
		(C3) If $z_1>z_0$, then $[z_0,z_1]\cap TC=\emptyset$. 
  
		Suppose $y^0$, $y^1$ are the optimal solutions of problems $(D^{z_0})$ and $(D^{z_1})$, respectively and $z^l_0$, $z^r_0$ are the left and right turning coordinates closest to $z_0$, respectively. 
		Let $k_j=\frac{y^0_j-y^1_j}{z_0-z_1}$ and $y^z_j=k_j(z-z_0)+y^0_j$, $j\in J\cup \bar J$. Notice that $(D^z_s)$ has the same optimal basis for any $z\in [z^l_0, z^r_0]$. Therefore, $y^z$ is the optimal solution of the problem $(D^z)$ for any $z\in [z^l_0, z^r_0]$ if and only if $y^z$ satisfies the constraints (\ref{D-z-con-2})-(\ref{D-z-con-6}), which can be formulated as $-d_j\leq y^z_j\leq 0$ for $j\in J$ and $x^0_jz-d_j\leq y^z_j\leq x^0_jz+d_j$ for $j\in \bar{J}$. Therefore, we have the following properties. $1)$ if $j\in J$ and $k_j>0$, then $\frac{k_jz_0-d_j-y^0_j}{k_j}\leq z\leq \frac{k_jz_0-y^0_j}{k_j}$. $2)$ if $j\in J$ and $k_j<0$, then $\frac{k_jz_0-y^0_j}{k_j}\leq z\leq \frac{k_jz_0-d_j-y^0_j}{k_j}$. $3)$ if $j\in \bar{J}$ and $k_j>x^0_j$, then $\frac{k_jz_0-d_j-y^0_j}{k_j-x^0_j}\leq z\leq \frac{k_jz_0+d_j-y^0_j}{k_j-x^0_j}$. $4)$ if $j\in \bar{J}$ and $k_j<x^0_j$, then $ \frac{k_jz_0+d_j-y^0_j}{k_j-x^0_j}\leq z\leq \frac{k_jz_0-d_j-y^0_j}{k_j-x^0_j}$. 
		Let $z^1_0=\max\limits_{j\in J,k_j>0}\frac{k_jz_0-d_j-y^0_j}{k_j}$, $z^2_0=\max\limits_{j\in J,k_j<0}\frac{k_jz_0-y^0_j}{k_j}$, $z^3_0=\max\limits_{j\in \bar{J},k_j>x^0_j}\frac{k_jz_0-d_j-y^0_j}{k_j-x^0_j}$, $z^4_0=\max\limits_{j\in \bar{J},k_j<x^0_j}\frac{k_jz_0+d_j-y^0_j}{k_j-x^0_j}$, $z^5_0=\min\limits_{j\in J,k_j>0}\frac{k_jz_0-y^0_j}{k_j}$, $z^6_0=\min\limits_{j\in J,k_j<0}\frac{k_jz_0-d_j-y^0_j}{k_j}$, $z^7_0=\min\limits_{j\in \bar{J},k_j>x^0_j}\frac{k_jz_0+d_j-y^0_j}{k_j-x^0_j}$ and $z^8_0=\min\limits_{j\in \bar{J},k_j<x^0_j}\frac{k_jz_0-d_j-y^0_j}{k_j-x^0_j}$. Hence, we have 
		\eq \label{z^l_0 and z^r_0}
		z^l_0=\max\{z^1_0, z^2_0, z^3_0, z^4_0\},\ \ \  
		z^r_0=\min\{z^5_0, z^6_0, z^7_0, z^8_0\}.
		\nq
		
		
		Based on the above analysis, we have the following lemma.
		\begin{lem}
			Let $k_j=\frac{y^0_j-y^1_j}{z_0-z_1}$ and $y^z_j=k_j(z-z_0)+y^0_j$, $j\in J\cup \bar J$. Then $y^z$ is the optimal solution of the problem $(D^z)$ for any $z\in [z^l_0, z^r_0]$, where $z^l_0, z^r_0$ are defined as in (\ref{z^l_0 and z^r_0}).
		\end{lem}
		
		\subsection{Calculate the left and right break points $z_l$ and $z_r$ }
		In this subsection, we present two algorithms \textbf{LBP} and \textbf{RBP} to calculate the  left and right break points $z_l$ and $z_r$ (if exist), respectively. 
		As they are similar, we only discuss the main idea to calculate $z_l$ in detail. 
		
		Firstly, we check the existence of the left break point $z_l$. Initialize $\underline{z}_l:=-(n+|\bar{J}|)d_{max}-1$ and $\bar{z}_l:=0$. If $\mathcal{F}_{D^{\underline{z}_l}}\neq \emptyset$, then the left break point $z_l$ dose not exist.  
		Otherwise, we can use a 
		binary search method to determine an interval $[\underline{z}_l,\bar{z}_l]$ which satisfies $\bar{z}_l-\underline{z}_l< \frac{\Delta}{2}$ and $\mathcal{F}_{D^{\underline{z}_l}}=\emptyset$, $\mathcal{F}_{D^{\bar{z}_l}}\neq \emptyset$. 
		Secondly, let $z_0=\bar{z}_l$ and $z_1=z_0+\frac{\Delta}{2}$, then we have $\mathcal{F}_{D^{{z}_1}}\neq \emptyset$ and $z_0$, $z_1$ satisfy the conditions (C1) and (C3) in subsection \ref{sub_section_cal_k_z}. 
		Finally, calculate the optimal solutions $y^{0}, y^{1}$ of the problems $(D^{z_0})$ and $(D^{z_1})$. Then we can get $z_l=z^l_0$ by (\ref{z^l_0 and z^r_0}).
		
		\begin{breakablealgorithm} \label{Alg-calculate-z_l}
			\caption{$z_l$=\textbf{LBP}($A,b,c,d,x^0,K$).}
				
				\begin{algorithmic}[1]
					\REQUIRE The coefficient matrix $A$, the vectors $b,c,d,x^0$ and a value $K$.
					\ENSURE  The left break point $z_l$.
					\STATE Let $\underline{z}_l:=-(n+|\bar{J}|)d_{max}-1$, $\bar{z}_l:=0$ and $\Delta:=\frac{1}{(2|\bar{J}|x^0_{max})^2}$.
					\IF  {$\mathcal{F}_{D^{\underline{z}_l}}\neq \emptyset$}
					\RETURN $z_l:=-\infty$.
					\ELSE
					\WHILE{$\bar{z}_l-\underline{z}_l \geq \frac{\Delta}{2}$}
					\STATE Let $z:=\frac{\underline{z}_l+\bar{z}_l}{2}$.
					\IF {$\mathcal{F}_{D^{z}}\neq \emptyset$}
					\STATE $\bar{z}_l:=z$.
					\ELSE
					\STATE $\underline{z}_l:=z$.
					\ENDIF
					\ENDWHILE
					\STATE Let $z_1:=\bar{z}_l+\frac{\Delta}{2}$, $z_0:=\bar{z}_l$ and $y^{0}, y^{1}$ be the optimal solutions of the problems $(D^{z_0})$ and $(D^{z_1})$, respectively.
				\RETURN $z_l:=z^l_0$ by (\ref{z^l_0 and z^r_0}).		
					\ENDIF 	\end{algorithmic}
		\end{breakablealgorithm}
		
		Next, we give the algorithm below to calculate $z_r$ similar to algorithm \ref{Alg-calculate-z_l}.
		\begin{breakablealgorithm} \label{Alg-calculate-z_r}
			\caption{$z_r$=\textbf{RBP}($A,b,c,d,x^0,K$).}
				
				\begin{algorithmic}[1]
					\REQUIRE The coefficient matrix $A$, the vectors $b,c,d,x^0$ and a value $K$.
					\ENSURE  The right break point $z_r$.
					\STATE Let $\bar{z}_r:=(n+|\bar{J}|)d_{max}+1$, $\underline{z}_r:=0$ and $\Delta:=\frac{1}{(2|\bar{J}|x^0_{max})^2}$.
					\IF  {$\mathcal{F}_{D^{\bar{z}_r}}\neq \emptyset$}
					\RETURN $z_r:=+\infty$.
					\ELSE
					\WHILE{$\bar{z}_r-\underline{z}_r \geq \frac{\Delta}{2}$}
					\STATE Let $z:=\frac{\underline{z}_r+ \bar{z}_r}{2}$.
					\IF {$\mathcal{F}_{D^{z}}\neq \emptyset$}
					\STATE $\underline{z}_r:=z$.
					\ELSE
					\STATE $ \bar{z}_r:=z$.
					\ENDIF
					\ENDWHILE
					\STATE Let $z_1:=\underline{z}_r-\frac{\Delta}{2}$, $z_0:=\underline{z}_r$ and $y^0, y^1$ be the optimal solutions of the problems $(D^{z_0})$ and $(D^{z_1})$, respectively.
				\RETURN $z_r:=z^r_0$ by (\ref{z^l_0 and z^r_0}).
					\ENDIF 	\end{algorithmic}
		\end{breakablealgorithm}

		\subsection{An algorithm to solve the  inverse problem (\textit{RIOVLP}$_1$)}
  
		In this subsection, we design an algorithm for the problem (\textit{RIOVLP}$_1$). 
		Based on the previous analysis, the problem (\textit{RIOVLP}$_1$) can be solved as long as the critical value $z^*$ is determined.
		
		Now we describe the main idea to calculated $z^*$. 
		Firstly, calculate the left and right break points $z_l,z_r$ by Algorithms \ref{Alg-calculate-z_l} and \ref{Alg-calculate-z_r}.  If $k^+_{z_l}\leq \delta$, then $z^*:=z_l$ is the critical value by Case (3.2) in Theorem \ref{th-z*}. Similarly, if $k^-_{z_r}\geq \delta$, then $z^*:=z_r$ is the critical value by Case (1.1) in Theorem \ref{th-z*}. Initialize $\bar{\tau}_0:=(n+|\bar{J}|)d_{max}+1$ and $\underline{\tau}_0=-\bar{\tau}_0$. If  $k_{\underline{\tau}_0}<\delta$ or  $k_{\bar{\tau}_0}>\delta$, then the problem (\textit{RIOVLP}$_1$) is infeasible. Now it comes to the case that the problem (\textit{RIOVLP}$_1$) is feasible and $z_l$ or $z_r$ is not the critical value, then we can determine an interval $[\underline{\tau}_\kappa, \bar{\tau}_\kappa]$ including $z^*$ by a binary search method which satisfies $k^-_{\bar{\tau}_\kappa}<\delta\leq k^+_{\underline{\tau}_\kappa}$. 
		In the $\kappa$-th iteration of the binary search method, if the current value $z_\kappa$ is a turning coordinate and $k^+_{z_\kappa}\leq\delta\leq k^-_{z_\kappa}$, 
		then $z^*=z_\kappa$ is the critical value by  Case (4.3) in Theorem \ref{th-z*}. Otherwise, the binary search method terminates when the length of  $|\bar{\tau}_\kappa-\underline{\tau}_\kappa|<\Delta$. In this case, there is only one turning coordinate in the interval $[\underline{\tau}_\kappa, \bar{\tau}_\kappa]$, which is just the critical value.
		
		Next, we discuss the following questions. (Q1) If $z\notin TC$ and $\mathcal{F}_{D^{z}}\neq \emptyset$, then how to calculate the slope $k_z$ of $\psi(z)$. 
		(Q2) If $z$ is the left break point, then how to calculate the right slope $k^+_z$  of $\psi(z)$. 
		(Q3) If $z$ is the right break point, then how to calculate the left slope $k^-_z$  of $\psi(z)$. 
		To answer the question (Q1), let $y^z$ be the optimal solution of the problem $(D^{z})$. 
		Calculate the optimal solution ($\bar\pi^z,\bar\alpha^z,\bar\beta^z$) of the problem  $(P^z)$ by duality theory of \textit{LP}. Hence, we have $k_z:=\sum_{j\in \bar{J}}(\bar\beta^z_j-\bar\alpha^z_j)x^0_j$.
		For the question (Q2), let $z':=z+\frac{\Delta}{2}$, then $z'\notin TC$ and $\mathcal{F}_{D^{z'}}\neq \emptyset$. Hence, we have $k^+_z:=k_{z'}$, where $k_{z'}$ can be calculated similar to the question (Q1). As for the question (Q3), let $z':=z-\frac{\Delta}{2}$, then $z'\notin TC$ and $\mathcal{F}_{D^{z'}}\neq \emptyset$. Hence, we have $k^-_z:=k_{z'}$, where $k_{z'}$ can be calculated similar to the question (Q1). For convenience, we give an algorithm to calculate $k_z$, $k^+_z$ or $k^-_z$.
		\begin{breakablealgorithm} \label{Alg-calculate-slope}
			\caption{$k_z$=\textbf{Slope}($z_l,z_r,z$).}
				
				\begin{algorithmic}[1]
					\REQUIRE The left and right break points $z_l$, $z_r$ and a value $z$.
					\ENSURE  The slope $k_z$ of $\psi(z)$.
					\IF{$z=z_l$}
					\STATE Let $z':=z+\frac{\Delta}{2}$.
					\ELSIF{$z=z_r$}
					\STATE Let $z':=z-\frac{\Delta}{2}$.
					\ELSE
					\STATE Let $z':=z$.
					\ENDIF
					\STATE Let $y^{z'}$ be the optimal solution of the problem $(D^{z'})$. Calculate the optimal solution ($\bar\pi^{z'},\bar\alpha^{z'},\bar\beta^{z'}$) of the problem  $(P^{z'})$ by duality theory of \textit{LP}.
					\RETURN $k_z:=\sum_{j\in \bar{J}}(\bar\beta^{z'}_j-\bar\alpha^{z'}_j)x^0_j$.
				\end{algorithmic}
		\end{breakablealgorithm}

		Next, we give Algorithm \ref{Alg-model-P} to solve the problem (\textit{RIOVLP}$_1$).
		\begin{breakablealgorithm} \label{Alg-model-P}
			\caption{$c^*$=\textbf{RIOVLP}($A,b,c,d,x^0,K$).}
				
				\begin{algorithmic}[1]
					\REQUIRE The coefficient matrix $A$, the vectors $b,c,d,x^0$ and a value $K$.
					\ENSURE  An optimal solution $c^*$ of the problem (\textit{RIOVLP}$_1$).
					\STATE Calculate $z_l:=$\textbf{LBP}($A,b,c,d,x^0,K$) and $z_r:=$\textbf{RBP}($A,b,c,d,x^0,K$).
					\STATE Initialize $\Delta:=\frac{1}{(2|\bar{J}|x^0_{max})^2}$, $\delta:=c_{\bar J}x^0_{\bar J}-K$,   $z^*:=+\infty$ and $\kappa:=0$.
					\IF{$z_l=-\infty$}
					\STATE Let $\underline{\tau}_0:=-(n+|\bar{J}|)d_{max}-1$ and $k_{\underline{\tau}_0}:=$\textbf{Slope}($z_l,z_r,\underline{\tau}_0$).
					\IF{$k_{\underline{\tau}_0}<\delta$}
					\STATE The problem (\textit{RIOVLP}$_1$) is infeasible and stop.
					\ENDIF
					\ELSE
					\STATE Let $\underline{\tau}_0:=z_l$ and $k^+_{\underline{\tau}_0}:=$\textbf{Slope}($z_l,z_r,\underline{\tau}_0$).
					\IF{$k^+_{\underline{\tau}_0}\leq \delta$}
					\STATE Let $z^*:=z_l$.
					\ENDIF
					\ENDIF
					
					\IF{$z_r=+\infty$}
					\STATE Let $\bar{\tau}_0:=(n+|\bar{J}|)d_{max}+1$ and $k_{\bar{\tau}_0}:=$\textbf{Slope}($z_l,z_r,\bar{\tau}_0$).
					\IF{$k_{\bar{\tau}_0}>\delta$}
					\STATE The problem (\textit{RIOVLP}$_1$) is infeasible and stop.
					\ENDIF
					\ELSE
					\STATE Let $\bar{\tau}_0:=z_r$ and $k^-_{\bar{\tau}_0}:=$\textbf{Slope}($z_l,z_r,\bar{\tau}_0$).
					\IF{$k^-_{\bar{\tau}_0}\geq \delta$}
					\STATE Let $z^*:=z_r$.
					\ENDIF
					\ENDIF
					\WHILE{$\bar{\tau}_{\kappa}-\underline{\tau}_{\kappa}\geq \Delta$ and $z^*=+\infty$}
					\STATE Let $z_{\kappa}:=\frac{\bar{\tau}_{\kappa}+\underline{\tau}_{\kappa}}{2}$, $z^1_{\kappa}:=z_{\kappa}-\frac{\Delta}{2}$, $z^2_{\kappa}:=z_{\kappa}-\frac{\Delta}{4}$, $z^3_{\kappa}:=z_{\kappa}+\frac{\Delta}{4}$ and  $z^4_{\kappa}:=z_{\kappa}+\frac{\Delta}{2}$.
					\STATE Let $k^1_{\kappa}:=\frac{\psi(z_{\kappa})-\psi(z^1_{\kappa})}{z_{\kappa}-z^1_{\kappa}}$, $k^2_{\kappa}:=\frac{\psi(z_{\kappa})-\psi(z^2_{\kappa})}{z_{\kappa}-z^2_{\kappa}}$, $k^3_{\kappa}:=\frac{\psi(z_{\kappa})-\psi(z^3_{\kappa})}{z_{\kappa}-z^3_{\kappa}}$ and $k^4_{\kappa}:=\frac{\psi(z_{\kappa})-\psi(z^4_{\kappa})}{z_{\kappa}-z^4_{\kappa}}$.
					\IF{$k^1_{\kappa}=k^2_{\kappa}$ and $k^3_{\kappa}=k^4_{\kappa}$ and $k^2_{\kappa}\neq k^3_{\kappa}$}
					\IF{$k^4_{\kappa}\geq \delta$}
					\STATE Update $\underline{\tau}_{\kappa+1}:=z^4_{\kappa}$ and $\bar{\tau}_{\kappa+1}:=\bar{\tau}_{\kappa}$.
					\ELSIF{$k^1_{\kappa}< \delta$}
					\STATE Update $\bar{\tau}_{\kappa+1}:=z^1_{\kappa}$ and $\underline{\tau}_{\kappa+1}:=\underline{\tau}_{\kappa}$.
					\ELSE
					\STATE Update $z^*:=z_{\kappa}$.
					\ENDIF
					\ELSE
					\STATE Let $k_{z_{\kappa}}:=$\textbf{Slope}($z_l,z_r,z_{\kappa}$).
					\IF{$k_{z_{\kappa}}\geq \delta$}
					\STATE Update $\underline{\tau}_{\kappa+1}:=z_{\kappa}$ and $\bar{\tau}_{\kappa+1}:=\bar{\tau}_{\kappa}$.
					\ELSE   
					\STATE Update $\bar{\tau}_{\kappa+1}:=z_{\kappa}$ and $\underline{\tau}_{\kappa+1}:=\underline{\tau}_{\kappa}$.
					\ENDIF
					\ENDIF
					\STATE Update $\kappa:=\kappa+1$.
					\ENDWHILE
					
					\IF {$z^*=+\infty$}
					\STATE Let $k_{\underline{\tau}_{\kappa}}:=$\textbf{Slope}($z_l,z_r,\underline{\tau}_{\kappa}$) and $k_{\bar{\tau}_{\kappa}}:=$\textbf{Slope}($z_l,z_r,\bar{\tau}^{\kappa}$).
					\STATE Calculate $z^*:=\frac{\psi(\bar{\tau}_{\kappa})-\psi(\underline{\tau}_{\kappa})-k_{\bar{\tau}_{\kappa}}\bar{\tau}_{\kappa}+k_{\underline{\tau}_{\kappa}}\underline{\tau}_{\kappa}}{k_{\underline{\tau}_{\kappa}}-k_{\bar{\tau}_{\kappa}}}$.
					\ENDIF
					\STATE Let $y^*$ be the optimal solution of the problem $(D^{z^*})$. Then $(y^*,z^*)$ be the optimal solution of the problem (\textit{DRIOVLP}$_1$). 
					\STATE Calculate the optimal solution ($\pi^*,\alpha^*,\beta^*$) of the problem  (\textit{RIOVLP}$^2_1$) by duality theory of \textit{LP}. 
					\RETURN $c^*:=c+\alpha^*-\beta^*$.

				\end{algorithmic}
		\end{breakablealgorithm}

		For convenience, we define $\mathcal{L}=\max\{d_{max},x^0_{max},n\}$. Then we can get the time complexity of Algorithm \ref{Alg-model-P}.
		\begin{thm} \label{Alg1-complexity}
			Algorithm \ref{Alg-model-P} can solve the problem (\textit{RIOVLP}$_1$) by solving the problem $(D^z)$ $O(\log \mathcal{L})$ times at most.
		\end{thm}
		\begin{proof} 
			The correctness of  Algorithm \ref{Alg-model-P} can be obtained by the main idea of the algorithm and Theorem \ref{th-z*}. Now we analyze the time complexity.  The main computation is in Line 1 and the while loop in Lines 25-45, which are all performed by a binary search method until the length of $|\bar{\tau}_\kappa-\underline{\tau}_\kappa|<\Delta$. We only need to calculate the number of iterations in the while loop. The initial interval length is $|\bar{\tau}_0-\underline{\tau}_0|=2(n+|\bar{J}|)d_{max}+2$ and the interval length will be reduced by at least 
			half in each iteration of a binary search method.
			Suppose there are $t$ iterations in the while loop in the worst case. Then we have $\big(2(n+|\bar{J}|)d_{max}+2\big)(\frac{1}{2})^t< \Delta=\frac{1}{(2|\bar{J}|x^0_{max})^2}$, which means $t> 2\log(2|\bar{J}|x^0_{max})+\log\big(2(n+|\bar{J}|)d_{max}+2\big)$.
			Hence, we have $t=\lceil 2\log(2|\bar{J}|x^0_{max})+\log\big(2(n+\bar{J})d_{max}+2\big)\rceil$  $\leq \lceil 2\log(2|\bar{J}|x^0_{max})+\log 3(n+\bar{J})d_{max}\rceil$=$O(\log \max\{d_{max},x^0_{max},n\}) = O(\log \mathcal{L})$.
			Notice that the Algorithm \ref{Alg-model-P} needs to calculate the problem $(D^z)$ at most five times in each iteration.
			Therefore, the conclusion holds. \hb
		\end{proof}

		\section{Applications to the Hitchcock and Shortest Path problems} \label{section_applications}
		In this section, we apply the previous research methods to the restricted inverse optimal value problems on Hitchcock and shortest path problem under weighted $l_1$ norm, respectively. 
  
  As these two problems can finally be transformed into a minimum cost flow  (\textit{MCF})  problem, we first introduce the problem  (\textit{MCF}) in \cite{Ahuja_book_1993}.
		
		Let $G(V,E,c,u)$ be a directed network with a cost $c_{ij}>0$ and a capacity $u_{ij}>0$ associated with every arc $(i,j)\in E$. We associate with each node $i\in V$ a supply $b(i)>0$ or a demand $b(i)<0$. Suppose $\sum_{i\in V}b(i)=0$, then 
the problem  (\textit{MCF}) can be stated as follows.
		\eqa  
		\min & & \sum_{(i,j)\in E}c_{ij}x_{ij} \nn\\
		 (\textit{MCF}) \ {\hbox {s.t.}} 
		&& \sum_{j:(i,j)\in E}x_{ij}-\sum_{j:(j,i)\in E}x_{ji}=b(i), i\in V, \nn\\
		&& 0\leq x_{ij}\leq u_{ij}, (i,j)\in E. \nn
		\nqa
		
		So far, the best strong polynomial time complexity for solving the problem  (\textit{MCF}) is $O(|E|\log |V|(|E|+|V|\log |V|))$ presented by \cite{Orlin_1993}.
		
		\subsection{The restricted inverse optimal value problem on Hitchcock problem under weighted $l_1$ norm}
		In this subsection, we study the restricted inverse optimal value problem (\textit{RIOVHC}$_1$) on Hitchcock problem  under weighted $l_1$ norm. 
		
		The Hitchcock problem can be described as follows. We have $m$ sources of some commodity, each with a supply of $a_i>0$ units, $i=1,\cdots,m$, and $n$ terminals, each with a demand of $b_j>0$ units, $j=1,\cdots,n$. Suppose $\sum_{i=1}^m a_i=\sum_{j=1}^n b_j$. There is a unit cost $c_{ij}>0$ of sending the commodity from source $i$ to terminal $j$. We aim to satisfy the demands at minimum cost. Hence, the Hitchcock problem (\textit{HC}) can be stated as follows.
		\eqa  
		\min & & \sum_{i=1}^m\sum_{j=1}^n c_{ij}x_{ij} \nn\\
		(\textit{HC}) \ {\hbox {s.t.}} 
		&& \sum_{j=1}^n x_{ij}=a_i, i=1,\cdots,m, \nn\\
		&& \sum_{i=1}^m x_{ij}=b_j, j=1,\cdots,n, \nn\\
		&& x_{ij} \geq 0. \nn
		\nqa

		Let $A$ and $c$ be the coefficient matrix and cost vector of the problem (\textit{HC}), respectively. The problem (\textit{RIOVHC}$_1$) can be described as follows. Given a feasible solution $x^0$ of the problem (\textit{HC}), a weight vector $d>0$ and a real number $K$,  we aim to adjust the cost vector $c$ to $\bar{c}$ under weighted $l_1$ norm such that $x^0$ becomes an optimal solution of the problem (\textit{HC}) under $\bar{c}$ and $\bar{c}x^0$ equals $K$.  
		
		Note that $A$ is a unimodular matrix. Suppose Assumption \ref{assumption-x^0} holds. 
		Let $J=\{(i,j)|x^0_{ij}=0\}$ and $\bar{J}=\{(i,j)|x^0_{ij}>0\}$. 
		Therefore, if we can solve the following problem (\textit{HC-}$\breve D^z$), then the problem (\textit{RIOVHC}$_1$) can also be solved by Algorithm \ref{Alg-model-P}.
		\eqa \label{HC-D-z}	
		\min && \sum_{(i,j)\in J\cup \bar{J}}c_{ij} x_{ij}-\sum_{(i,j)\in  \bar{J}}c_{ij} x^0_{ij}  \nn\\
		(HC\text{-} \breve D^z) \ {\hbox {s.t.}}  && \sum_{j=1}^n x_{ij}=a_i, i=1,\cdots,m, \nn\\
		&& \sum_{i=1}^m x_{ij}=b_j, j=1,\cdots,n, \nn\\  
		&&  0\leq x_{ij}\leq d_{ij}, (i,j)\in J,  \nn  \\
		&&  (1-z)x^0_{ij}-d_{ij}\leq x_{ij} \leq (1-z)x^0_{ij}+d_{ij}, (i,j)\in \bar{J}.   \nn
		\nqa
		
		Let $x'_{ij}=x_{ij}$ for any $(i,j)\in J$, $x'_{ij}=x_{ij}+d_{ij}-(1-z)x^0_{ij}$ for any $(i,j)\in \bar{J}$, $a'_i=a_i-\sum_{j:(i,j)\in \bar{J}} \big((1-z)x^0_{ij}-d_{ij} \big)$, $i=1,\dots,n$ and $b'_j=b_j-\sum_{i:(i,j)\in \bar{J}}\big((1-z)x^0_{ij}-d_{ij}\big)$, $j=1,\dots,n$. Hence, the problem (\textit{HC-}$\breve D^z$) can be turned into the following form.
		\eqa  	
		\min && \sum_{(i,j)\in J\cup \bar{J}}c_{ij} x'_{ij}-\sum_{(i,j)\in \bar{J}}c_{ij}(d_{ij}+zx^0_{ij}) \nn\\
		(HC\text{-}\breve{D}^z) \ {\hbox {s.t.}}  && \sum_{j:(i,j)\in J\cup \bar{J}} x'_{ij}=a'_i, i=1,\cdots,m, \nn\\
		&& \sum_{i:(i,j)\in J\cup \bar{J}}x'_{ij}=b'_j, j=1,\cdots,n, \nn\\  
		&&  0 \leq x'_{ij} \leq d_{ij}, (i,j)\in J,  \nn  \\
		&&  0 \leq x'_{ij} \leq 2d_{ij}, (i,j)\in \bar{J}.   \nn
		\nqa
		Obviously, we have $\sum_{i=1}^m a'_i=\sum_{j=1}^n b'_j$ since $\sum_{i=1}^m a_i=\sum_{j=1}^n b_j$. Furthermore, if there exists $a'_i<0$ or $b'_j<0$, then the problem (\textit{HC-}$\breve D^z$) is infeasible. Otherwise, the problem (\textit{HC-}$\breve D^z$) is a Hitchcock problem with upper bound constraints, which can be transformed into an (\textit{MCF}) problem. Therefore, we can obtain the time complexity of the problem (\textit{RIOVHC}$_1$). 
		
		\begin{thm}  
			The restricted inverse optimal value problem (\textit{RIOVHC}$_1$) on Hitchcock problem  under weighted $l_1$ norm  can be solved by Algorithm \ref{Alg-model-P} in $O\big((m\log n(m+n\log n))\log \mathcal{L}\big)$ time.
		\end{thm}
		
		Next, we present an example to execute Algorithm \ref{Alg-model-P} for the problem (\textit{RIOVHC}$_1$). 
		
\begin{figure}[htbp]  
			\centering
			\includegraphics[totalheight=6.5cm]{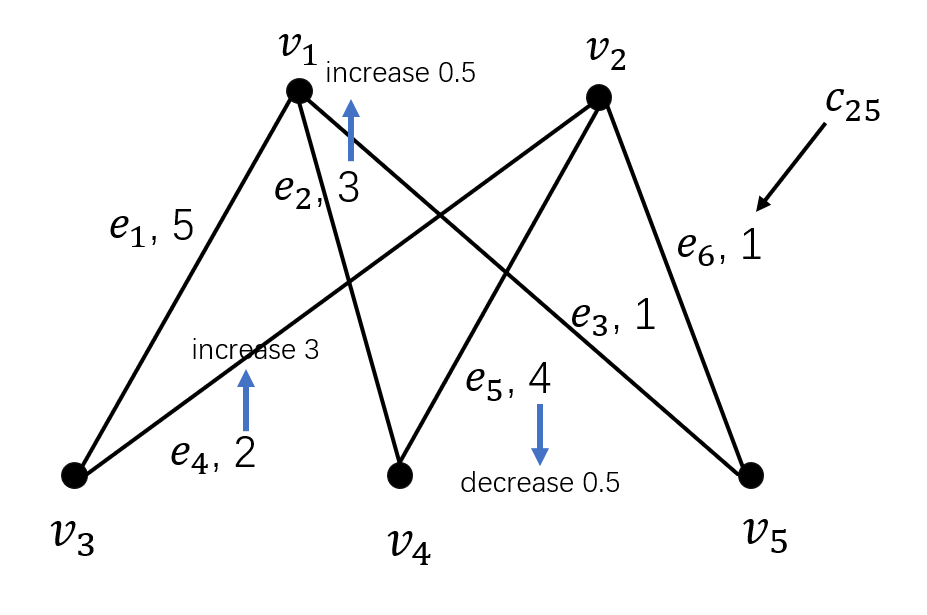}
			\caption{An example of the problem (\textit{RIOVHC}$_1$).}	\label{HC_graph}	
		\end{figure}

		\textbf{ Example 1.} 
		Let $v_1, v_2$ be the sources of some commodity with a supply of $a_1:=8, a_2:=12$ units, and three terminals $v_3$, $v_4$, $v_5$ with a demand of $b_1:=5, b_2:=4, b_3:=11$ units. There is a unit cost $c_{ij}>0$ of sending the commodity from source $i$ to terminal $j$ as shown in Figure \ref{HC_graph}. Let $x^0:=(3,2,3,2,2,8)$ be a feasible transportation strategy. Let $d:=(6,4,2,5,4,3)$ and $K:=50$. We aim to adjust the vector $c$ to $\bar{c}$ under weighted $l_1$ norm  such that $x^0$ becomes an optimal transportation strategy whose cost is just $K$ under $\bar{c}$.
		
\begin{figure}[htbp]  
			\centering
			\includegraphics[totalheight=7cm]{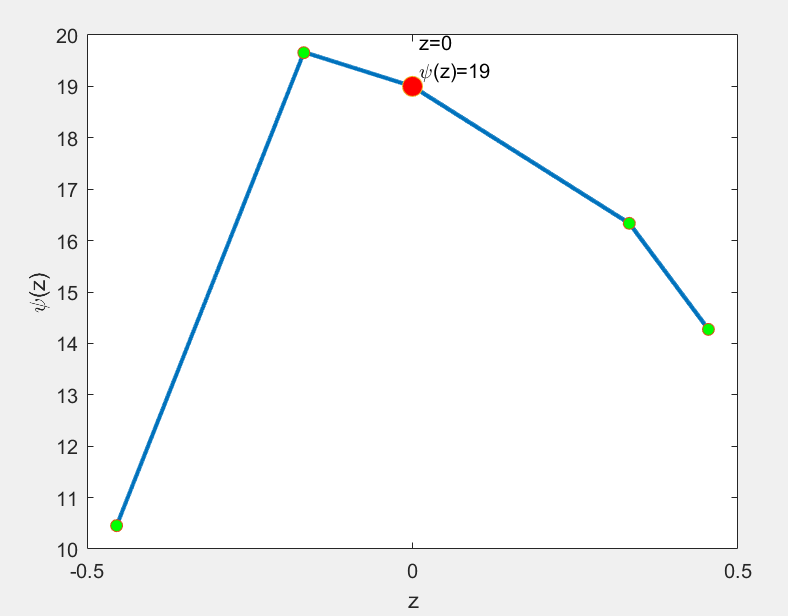}
			\caption{The graph of function $\psi(z)$ in Example 1.}	\label{HC_psi_z}	
		\end{figure}

		(1) Calculate $J:=\emptyset$, $\bar{J}:=\{1,2,3,4,5,6\}$, $d_{max}:=6$, $x^0_{max}:=8$, $\Delta:=\frac{1}{(2|\bar{J}|x^0_{max})^2}:=\frac{1}{9216}$ and $\delta:=c_{\bar J}x^0_{\bar J}-K:=-6$.
		(2) Calculate $z_l:=-\frac{5}{11}$ and $z_r:=\frac{5}{11}$ by Algorithms $\ref{Alg-calculate-z_l}$ and  $\ref{Alg-calculate-z_r}$.   
		(3) Let $\underline{\tau}_0:=z_l$ and $\bar{\tau}_0:=z_r$. Calculate $k^+_{\underline{\tau}_0}:=32$ and $k^-_{\bar{\tau}_0}:=-17$. 
		Therefore, the critical vale $z^*$ is located in the interval $[\underline{\tau}_0,\bar{\tau}_0]$ since $k^+_{\underline{\tau}_0}\geq \delta>k^-_{\bar{\tau}_0}$. We divide the interval $[\underline{\tau}_0,\bar{\tau}_0]$ by a binary search method. 
		Let $z_0:=\frac{\underline{\tau}_0+\bar{\tau}_0}{2}:=0$. Calculate $z^1_0:=-\frac{1}{18432}$, $z^2_0:=-\frac{1}{36864}$, $z^3_0:=\frac{1}{36864}$, $z^4_0:=\frac{1}{18432}$, and 
		$k^1_0:=k^2_0:=-4$, $k^3_0:=k^4_0:=-8$.  
		Hence, $z_0$ is a turning coordinate. 
		Notice that $k^4_0<\delta:=-6\leq k^1_0$. Therefore, $z_0$ is the critical value and the optimal solution of the problem (\textit{RIOVHC}$_1$) is $c^*:=(5,\textbf{3.5},1,\textbf{5},\textbf{3.5},1)$ as shown in Figure \ref{HC_graph}. Furthermore, we can draw the graph of function $\psi(z)$ by enumerating $z$ in the interval $[z_l,z_r]$ as shown in Figure \ref{HC_psi_z}, where the green dots represent the turning coordinates and the red dot is the critical value.
		
		\subsection{The restricted inverse optimal value problem on shortest path problem under weighted  $l_1$ norm}
		In this subsection, we study the restricted inverse optimal value problem on shortest path problem under weighted $l_1$ norm (\textit{RIOVSP}$_1$).
		
		Let $G=(V,E,c)$ be a directed network, where $V$, $E$ and $c$ denote the node set, the edge set and the edge cost vector, respectively.  
		Let nodes $s$ and $t$ denote two specified nodes. Suppose the network $G$ does not contain any negative cost cycle, then the $s-t$ shortest path problem can be described as follows. 
		\eqa  
		\min & & \sum_{(i,j)\in E}c_{ij}x_{ij} \nn\\
		(\textit{SP}) \ {\hbox {s.t.}} 
		&& \sum_{j:(i,j)\in E}x_{ij}-\sum_{j:(j,i)\in E}x_{ji}=1, i=s, \nn\\
		&& \sum_{j:(i,j)\in E}x_{ij}-\sum_{j:(j,i)\in E}x_{ji}=0, i \notin \{s,t\}, \nn\\
		&& \sum_{j:(i,j)\in E}x_{ij}-\sum_{j:(j,i)\in E}x_{ji}=-1, i=t, \nn\\
		&& x_{ij}\geq 0, (i,j)\in E. \nn
		\nqa
		
		The problem (\textit{RIOVSP}$_1$) can be described as follows.
		Let $P^0$ be a given $s-t$ path and $x^0$ be the corresponding 0-1 vector whose component 1 indicating the edges on $P^0$. Let $d> 0$ be the weight vector and $K$ be a real number. We aim to adjust the cost vector $c$ to $\bar{c}$ under weighted $l_1$ norm such that $x^0$ becomes a shortest path whose cost equals $K$ on new network $G=(V, E, \bar{c})$.   
		
		Obviously, the coefficient matrix of the problem (\textit{SP}) is unimodular. Suppose Assumption \ref{assumption-x^0} holds.
		Let $J=\{(i,j)|x^0_{ij}=0\}$ and $\bar{J}=\{(i,j)|x^0_{ij}=1\}$. 
		Then the problem (\textit{RIOVSP}$_1$) can also  be solved by Algorithm \ref{Alg-model-P} as long as the following problem  (\textit{SP-}$\breve D^z$) can be solved.
		\eqa  
		\min & & \sum_{(i,j)\in E}c_{ij}x_{ij}- \sum_{(i,j)\in   \bar{J}}c_{ij}\nn\\
		(SP\text{-}\breve D^z) \ {\hbox {s.t.}} 
		&& \sum_{j:(i,j)\in E}x_{ij}-\sum_{j:(j,i)\in E}x_{ji}=1, i=s, \nn\\
		&& \sum_{j:(i,j)\in E}x_{ij}-\sum_{j:(j,i)\in E}x_{ji}=0, i \notin \{s,t\}, \nn\\
		&& \sum_{j:(i,j)\in E}x_{ij}-\sum_{j:(j,i)\in E}x_{ji}=-1, i=t, \nn\\
		&& 0\leq x_{ij}\leq d_{ij}, (i,j)\in J,\nn\\
		&& 1-z-d_{ij}\leq x_{ij}\leq 1-z+d_{ij}, (i,j)\in \bar{J}. \nn
		\nqa
		
		Let $x'_{ij}=x_{ij}$ for any $(i,j)\in J$ and $x'_{ij}=x_{ij}-1+z+d_{ij}$ for any $(i,j)\in \bar{J}$. For convenience, we assume $P^0=j_0(s),j_1,\dots,j_k,j_{k+1}(t)$.  Let $b'(i)=1-(1-z-d_{ij_1})$ for $i=s$, $b'(i)=-(1-z-d_{ij_{h+1}})+(1-z-d_{j_{h-1}i})$ for $i=j_h$ and $1\leq h\leq k$, and $b'(i)=-1+(1-z-d_{j_ki})$ for $i=t$. 
		
		Hence, the problem (\textit{SP-}$\breve D^z$) can be turned into following form.
		\eqa \label{SP-D-z}	
		\min && \sum_{(i,j)\in E}c_{ij} x'_{ij}-\sum_{(i,j)\in \bar{J}}c_{ij}(z+d_{ij}) \nn\\
		(SP\text{-}\breve D^z) \ {\hbox {s.t.}}  && \sum_{j:(i,j)\in E} x'_{ij}-\sum_{j:(j,i)\in E} x'_{ji}=b'(i), i=s, \nn\\
		&& \sum_{j:(i,j)\in E} x'_{ij}-\sum_{j:(j,i)\in E} x'_{ji}=b'(i), i=j_h,h=1,\dots,k, \nn\\ 
		&& \sum_{j:(i,j)\in E} x'_{ij}-\sum_{j:(j,i)\in E} x'_{ji}=0, i\notin V(P^0), \nn\\
		&& \sum_{j:(i,j)\in E} x'_{ij}-\sum_{j:(j,i)\in E} x'_{ji}=b'(i), i=t, \nn\\
		&&  0 \leq x'_{ij} \leq d_{ij}, (i,j)\in J,  \nn  \\
		&&  0 \leq x'_{ij} \leq 2d_{ij}, (i,j)\in \bar{J}.   \nn
		\nqa
		
		Notice that for each edge $(j_{h},j_{h+1})\in P_0$ $(h=0,1,\cdots,k)$, there is one item  $-(1-z-d_{j_hj_{h+1}})$ for $i=j_h$ and one item $+(1-z-d_{j_{h}j_{h+1}})$ for $i=j_{h+1}$   in $b'(i)$. Then we have $\sum_{i\in V}b'(i)=\sum_{i\in V(P^0)}b'(i)=0$. Hence,  
		for a given $z\in \mathbb{R}$, the problem (\textit{SP-}$\breve D^z$) can be transformed into an  (\textit{MCF}) problem. Therefore, we can obtain the time complexity of the problem (\textit{RIOVSP}$_1$). 
		
		\begin{thm}  
			The restricted inverse optimal value problem (\textit{RIOVSP}$_1$) on shortest path under weighted $l_1$ norm  can be solved by Algorithm \ref{Alg-model-P} in  $O\big((m\log$ $ n(m+n\log n))\log \max\{d_{max},n\}\big)$ time. Furthermore, the time complexity can be reduced to  $O\big((m\log n(m+n\log n))\log n\big)$ under unit $l_1$ norm, where $d_{max}=1$.
		\end{thm}
		
		Next, we present an example to execute Algorithm \ref{Alg-model-P} for the problem  (\textit{RIOVSP}$_1$).
		
\begin{figure}[htbp]  
			\centering
			\includegraphics[totalheight=6.5cm]{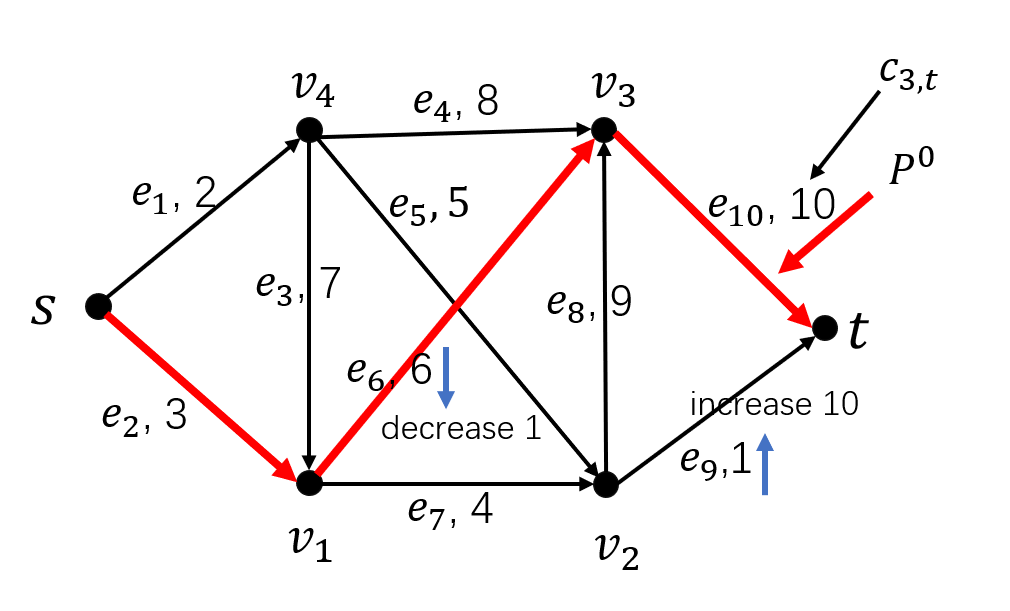}
			\caption{An example of the problem (\textit{RIOVSP}$_1$).}	\label{SP_graph}	
		\end{figure}

		\textbf{ Example 2.} 
		Let $G(V,E,c)$ be a directed weighted graph as shown in Figure \ref{SP_graph}, $P^0:=\{e_2,e_6,e_{10}\}$ (the read edges) be a given $s-t$ path, $c:=(2,3,7,8,5,6,4,9,1,10)$, $d:=(10,3,8,2,5,1,4,9,7,6)$ and $K:=18$. We aim to adjust the vector $c$ to $\bar{c}$ under weighted $l_1$ norm such that $P^0$ becomes a shortest $s-t$ path whose length is just $K$ under $\bar{c}$.
		
\begin{figure}[htbp]  
			\centering
			\includegraphics[totalheight=7cm]{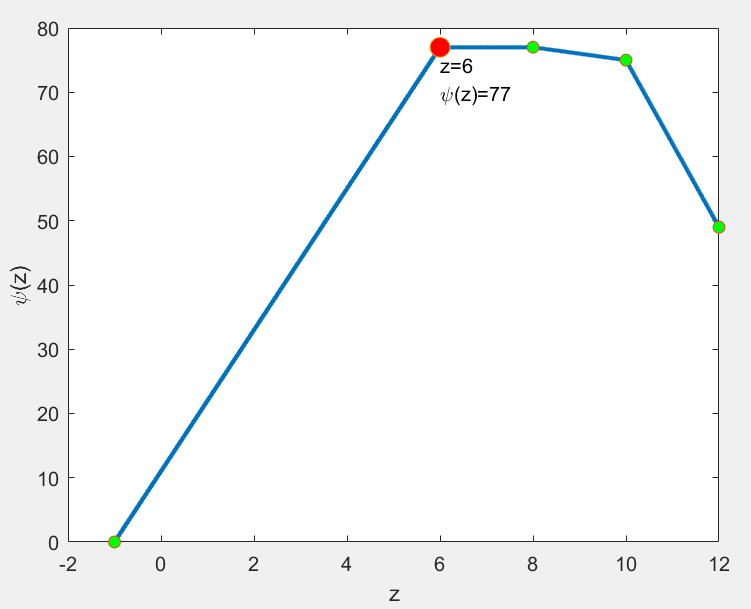}
			\caption{The graph of function $\psi(z)$ in Example 2.}	\label{SP_psi_z}	
		\end{figure}

		(1) Calculate $J:=\{1,3,4,5,7,8,9\}$, $\bar{J}:=\{2,6,10\}$, $d_{max}:=10$, $x^0_{max}:=1$, $\Delta:=\frac{1}{(2|\bar{J}|x^0_{max})^2}:=\frac{1}{36}$ and $\delta:=c_{\bar J}x^0_{\bar J}-K:=1$.
		(2) Calculate $z_l:=-1$ and $z_r:=12$ by Algorithms $\ref{Alg-calculate-z_l}$ and  $\ref{Alg-calculate-z_r}$.   
		(3) Let $\underline{\tau}_0:=z_l$ and $\bar{\tau}_0:=z_r$. Calculate $k^+_{\underline{\tau}_0}:=11$ and $k^-_{\bar{\tau}_0}:=-13$. 
		Therefore, the critical vale $z^*$ is located in the interval $[\underline{\tau}_0,\bar{\tau}_0]$ since $k^+_{\underline{\tau}_0}\geq \delta>k^-_{\bar{\tau}_0}$. We divide the interval $[\underline{\tau}_0,\bar{\tau}_0]$ by a binary search method. 
		Let $z_0:=\frac{\underline{\tau}_0+\bar{\tau}_0}{2}:=\frac{11}{2}$. Calculate $z^1_0:=\frac{395}{72}$, $z^2_0:=\frac{791}{144}$, $z^3_0:=\frac{793}{144}$, $z^4_0:=\frac{397}{72}$, and 
		$k^1_0:=k^2_0:=k^3_0:=k^4_0:=11$. 
		Hence, $z_0$ is not a turning coordinate. Calculate $k_{z_0}:=11$. Hence, $k_{z_0}\geq \delta$ and $\underline{\tau}_1:=z_0:=\frac{11}{2}$,  $\bar{\tau}_1:=\bar{\tau}_0:=12$. We continue to divide the interval $[\underline{\tau}_1,\bar{\tau}_1]$ by a binary search method. After nine iterations, we get the final interval $[\underline{\tau}_{9},\bar{\tau}_{9}]:=[\frac{3063}{512},\frac{769}{128}]$. (4) Calculate $k_{\bar{\tau}_{9}}:=0$, $k_{\underline{\tau}_{9}}:=11$, $\psi(\bar{\tau}_{9}):=77$ and $\psi(\underline{\tau}_{9}):=\frac{13902}{181}$. Hence, $z^*:=\frac{\psi(\bar{\tau}_{9})-\psi(\underline{\tau}_{9})-k_{\bar{\tau}_{9}}\bar{\tau}_{9}+k_{\underline{\tau}_{9}}\underline{\tau}_{9}}{k_{\underline{\tau}_{9}}-k_{\bar{\tau}_{9}}}:=6$. Therefore, the optimal solution of the problem (\textit{RIOVSP}$_1$) is $c^*:=(2,3,7,8,5,\textbf{5},4,9,\textbf{11},10)$ as shown in Figure \ref{SP_graph}. Furthermore, we can draw the graph of function $\psi(z)$ by enumerating $z$ in the interval $[z_l,z_r]$ as shown in Figure \ref{SP_psi_z}, where the green circles represent the turning coordinates and the red dot is the  critical value.
		
		\section{Conclusions and further research} \label{section_conclusions}
		In this paper, we mainly study the restricted inverse optimal value problem on (\textit{LP}) under weighted $l_1$ norm. Firstly, we construct the mathematical model of the problem 
  (\textit{RIOVLP}$_1$) by the dual theories, which is a linear programming problem. 
        Secondly, we introduce a sub-problem $(D^z)$  of the dual inverse problem (\textit{DRIOVLP}$_1$) with respect to a given value $z$ which only changes the upper and lower bounds of the variables compared to the original (\textit{LP}) problem. Thirdly, we design a binary search algorithm to calculate the critical value $z^*$ corresponding the optimal solution $(y^*,z^*)$ of the dual problem (\textit{DRIOVLP}$_1$). In each iteration, we need to solve a sub-problem problem $(D^z)$, which can be generally solved by an algorithm for the the original (\textit{LP}) problem. Finally, we can obtain an optimal solution of the inverse problem (\textit{RIOVLP}$_1$) by complementary slackness of \textit{LP}. The time complexity is $O(T^z\log \mathcal{L})$, where  $\mathcal{L}=\max\{d_{max},x^0_{max},n\}$ and $T^z$ is the time complexity to solve the sub-problem $(D^z)$. Finally, we apply the research methods to some restricted inverse optimal value problems on  Hitchcock and shortest path problems, where the sub-problem $(D^z)$ can be transformed into minimum cost flow problems. 
		
		We do not consider the bound constraints on  the adjustment amount in this paper, which may render some elements of the adjusted vector $\bar{c}$ too small or too large. In the future, we will study the bounded restricted inverse optimal value problem on \textit{LP} under weighted $l_1$ norm and other norms. Furthermore, if the original \textit{LP} problem is not standard, then we will consider whether our research results can be used to solve the corresponding inverse optimal value problem.
		
		
		
		
		\bibliographystyle{elsarticle-harv} 
		\bibliography{bibdatabase.bib}

\begin{thebibliography}{22}
\expandafter\ifx\csname natexlab\endcsname\relax\def\natexlab#1{#1}\fi
\providecommand{\url}[1]{\texttt{#1}}
\providecommand{\href}[2]{#2}
\providecommand{\path}[1]{#1}
\providecommand{\DOIprefix}{doi:}
\providecommand{\ArXivprefix}{arXiv:}
\providecommand{\URLprefix}{URL: }
\providecommand{\Pubmedprefix}{pmid:}
\providecommand{\doi}[1]{\href{http://dx.doi.org/#1}{\path{#1}}}
\providecommand{\Pubmed}[1]{\href{pmid:#1}{\path{#1}}}
\providecommand{\bibinfo}[2]{#2}
\ifx\xfnm\relax \def\xfnm[#1]{\unskip,\space#1}\fi
\bibitem[{Ahmed and Guan(2005)}]{Ahmed_2005}
\bibinfo{author}{Ahmed, S.}, \bibinfo{author}{Guan, Y.}, \bibinfo{year}{2005}.
\newblock \bibinfo{title}{The inverse optimal value problem}.
\newblock \bibinfo{journal}{Mathematical Programming} \bibinfo{volume}{102},
  \bibinfo{pages}{91--110}.
\bibitem[{Ahuja and Orlin(1993)}]{Ahuja_book_1993}
\bibinfo{author}{Ahuja, R.K.}, \bibinfo{author}{Orlin, J.B.},
  \bibinfo{year}{1993}.
\newblock \bibinfo{title}{Network Flows: Theory, Algorithms, and Applications}.
\newblock \bibinfo{publisher}{Prentice Hall}, \bibinfo{address}{New York}.
\bibitem[{Ahuja and Orlin(2001)}]{Ahuja_Orlin_2001}
\bibinfo{author}{Ahuja, R.K.}, \bibinfo{author}{Orlin, J.B.},
  \bibinfo{year}{2001}.
\newblock \bibinfo{title}{Inverse optimization}.
\newblock \bibinfo{journal}{Operations Research} \bibinfo{volume}{49},
  \bibinfo{pages}{771--783}.
\bibitem[{Burton and Toint(1992)}]{Burton_1992}
\bibinfo{author}{Burton, D.}, \bibinfo{author}{Toint, P.L.},
  \bibinfo{year}{1992}.
\newblock \bibinfo{title}{On an instance of the inverse shortest paths
  problem}.
\newblock \bibinfo{journal}{Mathematical Programming} \bibinfo{volume}{53},
  \bibinfo{pages}{45--61}.
\bibitem[{Chan and Kaw(2020)}]{Chan_kaw_2020}
\bibinfo{author}{Chan, T.C.}, \bibinfo{author}{Kaw, N.}, \bibinfo{year}{2020}.
\newblock \bibinfo{title}{Inverse optimization for the recovery of constraint
  parameters}.
\newblock \bibinfo{journal}{European Journal of Operational Research}
  \bibinfo{volume}{282}, \bibinfo{pages}{415--427}.
\bibitem[{Cui and Hochbaum(2010)}]{Cui_2010}
\bibinfo{author}{Cui, T.}, \bibinfo{author}{Hochbaum, D.S.},
  \bibinfo{year}{2010}.
\newblock \bibinfo{title}{Complexity of some inverse shortest path lengths
  problems}.
\newblock \bibinfo{journal}{Networks} \bibinfo{volume}{56},
  \bibinfo{pages}{20--29}.
\bibitem[{Ghobadi and Mahmoudzadeh(2021)}]{Ghobadi_Mahmoudzadeh_2021}
\bibinfo{author}{Ghobadi, K.}, \bibinfo{author}{Mahmoudzadeh, H.},
  \bibinfo{year}{2021}.
\newblock \bibinfo{title}{Inferring linear feasible regions using inverse
  optimization}.
\newblock \bibinfo{journal}{European Journal of Operational Research}
  \bibinfo{volume}{290}, \bibinfo{pages}{829--843}.
\bibitem[{Heuberger(2004)}]{Heuberger_2004}
\bibinfo{author}{Heuberger, C.}, \bibinfo{year}{2004}.
\newblock \bibinfo{title}{Inverse combinatorial optimization: A survey on
  problems, methods, and results}.
\newblock \bibinfo{journal}{Journal of Combinatorial Optimization}
  \bibinfo{volume}{8}, \bibinfo{pages}{329--361}.
\bibitem[{Huang and Liu(1999)}]{Huang_Liu_1999}
\bibinfo{author}{Huang, S.}, \bibinfo{author}{Liu, Z.}, \bibinfo{year}{1999}.
\newblock \bibinfo{title}{On the inverse problem of linear programming and its
  application to minimum weight perfect k-matching}.
\newblock \bibinfo{journal}{European Journal of Operational Research}
  \bibinfo{volume}{112}, \bibinfo{pages}{421--426}.
\bibitem[{Jia et~al.(2023)Jia, Guan, Wang, Zhang and Pardalos}]{Jia_Guan_2023}
\bibinfo{author}{Jia, J.}, \bibinfo{author}{Guan, X.}, \bibinfo{author}{Wang,
  H.}, \bibinfo{author}{Zhang, B.}, \bibinfo{author}{Pardalos, P.M.},
  \bibinfo{year}{2023}.
\newblock \bibinfo{title}{Combinatorial algorithms for solving the restricted
  bounded inverse optimal value problem on minimum spanning tree under weighted
  $l_\infty$ norm}.
\newblock \bibinfo{journal}{Journal of Computational and Applied Mathematics}
  \bibinfo{volume}{419}, \bibinfo{pages}{114754}.
\bibitem[{Lv et~al.(2010)Lv, Chen and Wan}]{LV_2010}
\bibinfo{author}{Lv, Y.}, \bibinfo{author}{Chen, Z.}, \bibinfo{author}{Wan,
  Z.}, \bibinfo{year}{2010}.
\newblock \bibinfo{title}{A penalty function method based on bilevel
  programming for solving inverse optimal value problems}.
\newblock \bibinfo{journal}{Applied Mathematics Letters} \bibinfo{volume}{23},
  \bibinfo{pages}{170--175}.
\bibitem[{Lv et~al.(2008)Lv, Hu and Wan}]{LV_2008}
\bibinfo{author}{Lv, Y.}, \bibinfo{author}{Hu, T.}, \bibinfo{author}{Wan, Z.},
  \bibinfo{year}{2008}.
\newblock \bibinfo{title}{A penalty function method for solving inverse optimal
  value problem}.
\newblock \bibinfo{journal}{Journal of Computational and Applied Mathematics}
  \bibinfo{volume}{220}, \bibinfo{pages}{175--180}.
\bibitem[{Mohajerin~Esfahani et~al.(2018)Mohajerin~Esfahani,
  Shafieezadeh-Abadeh, Hanasusanto and Kuhn}]{Mohajerin_2018}
\bibinfo{author}{Mohajerin~Esfahani, P.}, \bibinfo{author}{Shafieezadeh-Abadeh,
  S.}, \bibinfo{author}{Hanasusanto, G.A.}, \bibinfo{author}{Kuhn, D.},
  \bibinfo{year}{2018}.
\newblock \bibinfo{title}{Data-driven inverse optimization with imperfect
  information}.
\newblock \bibinfo{journal}{Mathematical Programming} \bibinfo{volume}{167},
  \bibinfo{pages}{191--234}.
\bibitem[{Orlin(1993)}]{Orlin_1993}
\bibinfo{author}{Orlin, J.B.}, \bibinfo{year}{1993}.
\newblock \bibinfo{title}{A faster strongly polynomial minimum cost flow
  algorithm}.
\newblock \bibinfo{journal}{Operations Research} \bibinfo{volume}{41},
  \bibinfo{pages}{338--350}.
\bibitem[{Papadimitriou and Steiglitz(1998)}]{Papadimitriou_book_1982}
\bibinfo{author}{Papadimitriou, C.}, \bibinfo{author}{Steiglitz, K.},
  \bibinfo{year}{1998}.
\newblock \bibinfo{title}{Combinatorial Optimization: Algorithms and
  Complexity}.
\newblock \bibinfo{publisher}{Dover Publications Inc}, \bibinfo{address}{New
  York}.
\bibitem[{Wang et~al.(2021)Wang, Guan, Zhang and Zhang}]{WangHui_2021}
\bibinfo{author}{Wang, H.}, \bibinfo{author}{Guan, X.}, \bibinfo{author}{Zhang,
  Q.}, \bibinfo{author}{Zhang, B.}, \bibinfo{year}{2021}.
\newblock \bibinfo{title}{Capacitated inverse optimal value problem on minimum
  spanning tree under bottleneck hamming distance}.
\newblock \bibinfo{journal}{Journal of Combinatorial Optimization}
  \bibinfo{volume}{41}, \bibinfo{pages}{861--887}.
\bibitem[{Zhang et~al.(2021)Zhang, Guan, Pardalos, Wang, Zhang, Liu and
  Chen}]{Zhang_Guan_2021}
\bibinfo{author}{Zhang, B.}, \bibinfo{author}{Guan, X.},
  \bibinfo{author}{Pardalos, P.M.}, \bibinfo{author}{Wang, H.},
  \bibinfo{author}{Zhang, Q.}, \bibinfo{author}{Liu, Y.},
  \bibinfo{author}{Chen, S.}, \bibinfo{year}{2021}.
\newblock \bibinfo{title}{The lower bounded inverse optimal value problem on
  minimum spanning tree under unit $l_\infty$ norm}.
\newblock \bibinfo{journal}{Journal of Global Optimization}
  \bibinfo{volume}{79}, \bibinfo{pages}{757--777}.
\bibitem[{Zhang et~al.(2020)Zhang, Guan and Zhang}]{Zhang_Guan_2020}
\bibinfo{author}{Zhang, B.}, \bibinfo{author}{Guan, X.},
  \bibinfo{author}{Zhang, Q.}, \bibinfo{year}{2020}.
\newblock \bibinfo{title}{Inverse optimal value problem on minimum spanning
  tree under unit $l_\infty$ norm}.
\newblock \bibinfo{journal}{Optimization Letters} \bibinfo{volume}{14},
  \bibinfo{pages}{2301--2322}.
\bibitem[{Zhang and Cai(1998)}]{Zhang_Cai_1998}
\bibinfo{author}{Zhang, J.}, \bibinfo{author}{Cai, M.}, \bibinfo{year}{1998}.
\newblock \bibinfo{title}{Inverse problem of minimum cuts}.
\newblock \bibinfo{journal}{Mathematical Methods of Operations Research}
  \bibinfo{volume}{47}, \bibinfo{pages}{51--58}.
\bibitem[{Zhang and Liu(1996)}]{Zhang_Liu_1996}
\bibinfo{author}{Zhang, J.}, \bibinfo{author}{Liu, Z.}, \bibinfo{year}{1996}.
\newblock \bibinfo{title}{Calculating some inverse linear programming
  problems}.
\newblock \bibinfo{journal}{Journal of Computational and Applied Mathematics}
  \bibinfo{volume}{72}, \bibinfo{pages}{261--273}.
\bibitem[{Zhang and Liu(1999)}]{Zhang_Liu_1999}
\bibinfo{author}{Zhang, J.}, \bibinfo{author}{Liu, Z.}, \bibinfo{year}{1999}.
\newblock \bibinfo{title}{A further study on inverse linear programming
  problems}.
\newblock \bibinfo{journal}{Journal of Computational and Applied Mathematics}
  \bibinfo{volume}{106}, \bibinfo{pages}{345--359}.
\bibitem[{Zhang et~al.(2023)Zhang, Guan, Jia, Qian and Pardalos}]{ZhangQ_2023}
\bibinfo{author}{Zhang, Q.}, \bibinfo{author}{Guan, X.}, \bibinfo{author}{Jia,
  J.}, \bibinfo{author}{Qian, X.}, \bibinfo{author}{Pardalos, P.M.},
  \bibinfo{year}{2023}.
\newblock \bibinfo{title}{The restricted inverse optimal value problem on
  shortest path under $l_1$ norm on trees}.
\newblock \bibinfo{journal}{Journal of Global Optimization}
  \bibinfo{volume}{86}, \bibinfo{pages}{251--284}.

\end{thebibliography}
		

		

	\end{document}